\documentstyle[11pt]{article}
\newtheorem{defn}{Definition}[section]
\newtheorem{lem}[defn]{Lemma}
\newtheorem{thm}[defn]{Theorem}
\newtheorem{prop}[defn]{Proposition}
\newtheorem{cor}[defn]{Corollary}
\newtheorem{rem}[defn]{Remark}

\newtheorem{assu}[defn]{Assumption}

\newtheorem{nota}[defn]{Notation}
\newtheorem{nota-rem}[defn]{Notation--Remark}
\newcommand{\C}{{\bf C}}
\newcommand{\f}{{\bf F}_3}
\newcommand{\epsi}{\epsilon}
\newcommand{\E}{{\cal E}}

\newcommand{\Ho}{H^0}

\newcommand{\im}{\mbox{Im}}

\newcommand{\k}{K_S}
\newcommand{\hS}{\hat{S}}
\newcommand{\hC}{\hat{C}}
\newcommand{\mod}{\,{\rm mod}\,}

\newcommand{\spec}{{\rm spec}}
\newcommand{\pp}{{\bf P}}
\newcommand{\si}{\sigma_{\infty}}
\newcommand{\so}{\sigma_0}
\newcommand{\qed}{$\quad \diamond$\par\smallskip}
\newcommand{\proof}{{\bf Proof:}\,\,}
\newcommand{\inv}{^{-1}}
\newcommand{\OO}{{\cal O}}
\newcommand{\OS}{\OO_S}
\newcommand{\OSh}{\OO_{\hS}}
\newcommand{\OC}{{\OO_C}}
\newcommand{\OX}{{\OO_X}}
\newcommand{\Of}{\OO_{\f}}
\newcommand{\proj}{\hbox{Proj}}

\title{Triple canonical surfaces of minimal degree}
\author{Margarida Mendes Lopes -- Rita Pardini}
\date{}
\begin{document}
\def\theequation{\thesection.\arabic{equation}}
\maketitle

\setcounter{defn}{0} \setcounter{equation}{0}

\section{Introduction} In this paper, continuing the study
of pathologies of the canonical map for surfaces of
general type that we started in \cite{mp}, we study the surfaces with
$p_g\geq 4$ for
which the canonical map is $3$--$1$ onto a surface of minimal
degree $n-1$ in $\pp^{n}$, $n\ge 3$, under the assumption that the
canonical system contains a smooth curve.

 We obtain a complete classification and we construct examples of all the
possible numerical cases. Contrarily to what happens for the similar problem
for $2$--$1$ canonical map, where there are examples with unbounded invariants
(for instance surfaces with $K^2=2p_g-4$, see \cite{ho-1}, \cite{pe}),
 it turns out that these surfaces only exist in a very limited range, namely
for $p_g\leq 5$ and $K^2\leq 9$.

This problem has been studied classically, in a slightly
different context, by Pompilj \cite{po}, for the case of
canonical system base point free.

More recently Starnone in his thesis, \cite{sta}, classified the surfaces
with $p_g\ge
4$ such
that the canonical system has at most simple distinct base points and the
general
canonical
curve is trigonal: it turns out that these surfaces are mapped $3$--$1$ onto a
surface  of minimal degree in
$\pp^{p_g-1}$ and therefore satisfy our assumptions.
 Our results, which were obtained at the same time and
independently from his, are more general, since we only
assume the existence of a smooth canonical curve,  and indeed there exist
surfaces
satisfying our assumptions but not Starnone's, as it is shown by our main
result:
\begin{thm}\label{Main} a) Let $S$ be a minimal surface of
general type such that
$p_g\ge 4$, the general canonical curve is smooth, and
the canonical map of $S$ is $3$--$1$ onto a surface
$\Sigma$ of minimal degree in
$\pp^{p_g-1}$. Then   $S$ is regular, $\Sigma$ is a cone,
$|\k|$ has at most $2$ base points; the invariants of $S$
satisfy
$$\matrix{ && \k^2 &&  p_g\cr (M_1)&&9&&5\cr (M_2)&&8&&4\cr
(M_3)&&7&&4\cr (M_4)&&6&&4\cr (N)&&8&&4\cr}.$$
 Furthermore $S$ contains a rational pencil $|C|$ such that the general
curve of $|C|$ is smooth, nonhyperelliptic and such that,
in cases $(M)$, $C^2=1$, $\k C=3$ and, in case $(N)$,
$C^2=2$, $\k C=4$.

In case $(N)$, $\k=2C$, whilst in case $(M_2)$, which has
the same numerical invariants, $\k$ is not 2-divisible in
$Pic (S)$. Furthermore in case $(N)$, $|\k|$ has two
simple base points which are infinitely near.

b) Conversely, any minimal surface $S$ of general type
with $p_g\geq 4$, containing a rational pencil $|C|$ such
that the general curve of $|C|$ is smooth, nonhyperelliptic, with $C^2=1$, $\k
C=3$ satisfies the hypothesis in a) and is one the surfaces of type $(M)$. On
the other hand if the minimal surface $S$ of general type, with $p_g\geq 4$,
contains a
rational pencil $|C|$ such that the general curve of $|C|$
is smooth, nonhyperelliptic, with $C^2=2$, $\k C=4$,
then necessarily $\k^2=8$, $\k=2C$, $p_g=4$, and the canonical map of $S$ is
generically finite of degree  $3$ or $4$. If
the degree is
$3$, then
$S$ satisfies the hypothesis in a) and is a surface of type $(N)$.

c) All the surfaces described in a) are generically finite
triple covers of the cubic cone in $\pp^4$, and it is
possible to construct examples of all such surfaces.

\end{thm}

The surfaces of types $(M_1)$  and $(M_4)$ were known to
Pompilj, who described them in terms of triple covers of
$\pp^2$, (see \cite{po}). More recently, Horikawa, doing
his classification of surfaces with $p_g=4$, $K^2=6$,
described and proved the existence of  surfaces of type
$(M_4)$ (see surfaces of type II in \cite{ho-3}), whilst
K. Konno constructed surfaces of type $(M_1)$ in
\cite{Konno}. Also  Zucconi, in his thesis,\cite{zu}, constructs
examples of surfaces of types $(M_3)$ and $(M_1)$.

Starnone, in his thesis, working under the assumption that
$|\k|$ has no infinitely near points, finds and constructs
all the surfaces of type  $(M)$. His results miss the type
$(N)$ surfaces exactly because for these surfaces $|\k|$
 always has infinitely near base points.  Although both
 Starnone's constructions  for type $(M)$ surfaces and ours
are done via triple covers, they differ in that Starnone
uses exactly the triple covers obtained from the canonical
map, whilst we use triple covers of the cubic cone in $\pp^3$ induced by
the linear
system $|3C|$, where $|C|$ is the pencil in thm. \ref{Main}.

We do not know any reference in the literature to surfaces
of type $(N)$.

The paper is organized as follows:

In section $2$ we show that that every surface $S$ with $3$--$1$
canonical map satisfying our assumptions contains a pencil
of curves as in \ref{Main} a).

In sections $3$ and $4$ we
study the surfaces having such a pencil and establish
various properties of these surfaces. In particular we show that these
surfaces are triple covers of the cubic cone
in $\pp^4$.

In section $5$ we describe these triple covers precisely, and finally in
section $6$ we give explicit examples of all these surfaces.

{\em Notation and conventions:} All varieties are  projective varieties
over the complex numbers.
The $n$-dimensional projective space is denoted by $\pp^n$.  As usual,
${\cal O}_Y$ is the structure sheaf of the variety
$Y$,
$H^i(Y,{\cal F})$ is the $i$-th cohomology group of a sheaf ${\cal
F}$ on $Y$, and $h^i(Y,{\cal F})$ is the dimension of $H^i(Y,{\cal
F})$; for a line bundle $M$ on $Y$, we denote by $|M|$ the complete
linear system $\pp(H^0(Y,M))$. When dealing with smooth varieties, we
do not distinguish between line bundles and divisors. If
$S$ is a smooth surface, then
$K_S$ denotes a canonical divisor,
$p_g(S)=h^0(S,\OS(K_S))$ is the {\em geometric genus} and
$q(S)=h^1(S, {\cal O}_S)$ is the {\em irregularity}, $K_S^2$ is the
self-intersection of the canonical divisor. A surface
$S$ is said to be {\em regular} if $q(S)= 0$. A (rational or linear) pencil
on a surface $S$ is a linear system on $S$ of projective dimension $1$.
The intersection
number of two divisors
$C$,
$D$ on a smooth surface is denoted simply by $CD$, linear equivalence is
denoted by $\equiv$ and numerical equivalence is denoted by $\sim$. We will
 sometimes refer to an effective non zero divisor on a
surface as to a ``curve''.
\bigskip

{\em Acknowledgements:} We wish to thank Patrizia Gianni and  Barry Trager,
who have shown great patience in teaching us  how to do computations with
Axiom, and have thus enabled us to construct the examples  of section $6$.

The present collaboration takes place in the framework of the HCM contract
AGE, no. ERBCHRXCT940557.  The first author is a member of CAUL and
of the project Praxis 2/2.1/MAT/73/94. The second author is a member
of GNSAGA of CNR. This paper was started during a visit, financed  by CAUL, of
the second author to Lisbon.

\setcounter{defn}{0} \setcounter{equation}{0}

\section{Surfaces whose canonical map is $3$--to--$1$ onto
a surface of minimal degree}\label{3-1}
\setcounter{defn}{0}
\setcounter{equation}{0}

In this section we prove the following:

\begin{thm}\label{III} Let $S$ be a minimal surface of
general type with
$p_g(S)\ge 4$, let
$\phi:S\to\pp^{p_g-1}$ be the canonical map of $S$, and let
$\Sigma$ be the image of $S$ via $\phi$. Suppose that the
general curve of
$|K_S|$ is smooth, that $\Sigma$ is a surface, that
$\deg\phi$ is equal to $3$ and that the degree of
$\Sigma$ is equal to
$p_g-2$, namely it is the lowest possible. Then either $6\le
\k^2\le 8$ and $p_g(S)=4$ or $\k^2=9$ and $p_g(S)=5$.  Furthermore $S$
contains a
rational pencil $|C|$ such that the general curve of
$|C|$ is irreducible and non-hyperelliptic, and such that
either $C^2=1$,
$\k C=3$ or $C^2=2$, $\k C=4$.
\end{thm}

\begin{rem} The requirement that the general canonical
curve be smooth is equivalent to the fact that $|K_S|$ has
no fixed component and that the base points (if any) are
simple,  possibly  infinitely near.
\end{rem}
\noindent We shall use the following general fact:

\begin{lem}\label{francia}  Let $X$ be a surface, let
$D$ be a curve on $X$ and let  $x_1,\ldots x_d\in {\rm
Sing}\,D$ be distinct points. Let $p: \tilde S \to S$ be
the blow-up at
$x_1,\ldots x_d$ and let $ E_1,\ldots E_d$ be the
corresponding exceptional curves. Setting
$D^\prime=p^*D-E_1-\cdots-E_d$ and
$D^\prime{}^\prime=p^*D-2E_1-\cdots-2E_d$, the following
two conditions are equivalent:

(i)  $x_1,\ldots x_d$ do not impose independent conditions
on  $|K_X + D|$.

(ii) The restriction map
$H^0(D^\prime,{\cal O}_{D^\prime})\to
H^0(D^\prime{}^\prime,{\cal O}_{D^\prime{}^\prime})$ is
not surjective.
\end{lem}
 \proof  The proof of this lemma is exactly like the
proofs of the analogous statements for $d=1$ in
\cite{francia} and $d=2$ in \cite{ml} and so we omit it.
\qed

 \noindent{\bf Proof of Theorem
\ref{III}:}  The proof consists of several steps.

\noindent {\em Claim $1$: $6\le \k^2\le 9$ and $p_g=4,5$.}

\noindent  Let $z\in\Sigma$ be a general point, let
$\phi\inv(z)=\{x_1, x_2,x_3\}$ and let $C$ be the
pull--back on $S$ of a general hyperplane section $H$ of
$\Sigma$ through $z$. The curve $H$ is  smooth rational
by the assumption that
$\Sigma$ has minimal degree. Therefore, since
$C$ is a smooth canonical curve, the linear system
$|x_1+x_2+x_3|$ on $C$ has dimension $1$. By Riemann--Roch
on $C$, the points $x_1, x_2,x_3$ fail to impose
independent conditions on $K_C$, and therefore, a
fortiori, also on
$|2K_S|$. Let
$H'$ be a tangent hyperplane to $\Sigma$ at $z$ and let
$D$ be the pull--back of
$H'$ to
$\Sigma$: the curve
$D$ is singular at $x_1, x_2,x_3$. By lemma \ref{francia},
using the notation introduced there,
$h^0(D'',\OO_{D''})\ge 2$ and thus, in particular $D''$ is
not 1--connected. Arguing exactly  as in the proof of thm.
$5$ of
\cite{bombieri}, one shows that there exist curves $A$ and
$B$ such that
$D=A+B$,
$AB= 2$ and at least two of the $x_i$'s belong to $A\cap
B$. \par
 Since $z$ can be chosen to be sufficiently general, one
can argue as in the
 proof of thm. 5 of
\cite{bombieri} and show that $\k A, \k B\ge 2$ and that
$\k A=2$ or $\k B=2$ iff
$S$ has a pencil of genus $2$. The latter possibility is
excluded,  since  the degree of the canonical map is
$3$.  So we can assume that
$\k A\ge 3$,
$\k B\ge 3$.

Say $\k A\le \k B$.  By the index theorem, $A^2B^2-(AB)^2\le
0$, with equality holding if and only if $B\sim mA$ for
some
$m\in Q$. Since
$\k A=A^2+AB=A^2+2$ and
$\k B=B^2+AB=B^2+2$ and
$\k$ is nef,  the numerical possibilities are
$$\matrix{ && \k^2 &&  \k  A && A^2 && B^2\cr
(M_1)&&9&&3&&1&&4\cr (M_2)&&8&&3&&1&&3\cr
(M_3)&&7&&3&&1&&2\cr (M_4)&&6&&3&&1&&1\cr
(N)&&8&&4&&2&&2\cr}.$$ where in case $(M_1)$, $B\sim 3A$
and in case $(N)$
$B\sim A$.

So $\k^2\le 9$. Moreover one has
\begin{equation}\label{inequality}
 K^2_S\ge \deg\phi \deg\Sigma=3(p_g(S)-2)
\end{equation} with equality holding if and only if
$|K_S|$ is base point free. Since $p_g(S)\ge 4$ by
assumption, one obtains immediately the assertion about
$p_g(S)$.
\medskip

\noindent {\em Claim
$2$: if
$p_g=5$, then
$\k^2=9$,
$\Sigma$ is the rational cubic cone and $S$ contains a
rational pencil $|C|$ such that the general curve of
$|C|$ is irreducible and non-hyperelliptic with $C^2=1$,
$\k C=3$. }

\noindent If $p_g=5$, by Claim $1$ and \ref{inequality} we
have $\k^2=9$ and
$|K_S|$ base point free. We wish to show that $\Sigma$ is
a cone. Assume on the contrary that
$\Sigma$ is a rational normal cubic scroll and denote by
$C$ the pull--back of a ruling of
$\Sigma$: one has $C^2=0$ and $\k C=3$, contradicting the
adjunction formula. So
$\Sigma$ is the rational normal cubic cone; denote again
by $C$  the pull--back of a ruling of $\Sigma$ and
consider the pull-back $D$ of a hyperplane section of
$\Sigma$ passing through  the vertex. One can write
$D=3C+Z$, where $Z$ is an effective divisor, possibly
empty, that is contracted to the vertex of the cone. Now
$9=\k^2=\k D=3\k C+\k Z=9+\k Z$, and so $\k Z=0$ and $\k
C=3$. Since $|C|$ is positive dimensional and the general $C$ is irreducible,
$C^2\ge 0$, and $C^2\ne 0$ by parity. So the Hodge index
theorem finally yields
$C^2=1$ and $\k\sim 3C$. Since $|\k|$ cuts on the general
curve in $|C|$ a $g_3^1$ without base points, this curve
 is not hyperelliptic.
\medskip

\noindent {\em Claim $3$: if $p_g=4$, then either  $S$
  contains a rational pencil $|C|$  with $C^2=1$,
$\k C=3$,  such that the general
curve of
$|C|$ is irreducible and non-hyperelliptic, or $\k^2=8$ and  $S$ contains a
rational
pencil
$|C|$ with $C^2=2$,
$\k C=4$, such that the general curve of
$|C|$ is irreducible and non-hyperelliptic. }

\noindent We start by showing that $\Sigma$ is the quadric
cone. Assume on the contrary that
$\Sigma$ is a smooth quadric, and denote by $D$ the
pull--back of a general tangent section of
$\Sigma$.  One can write
$D=A+B$, where
$A$ and
$B$ are pull--backs of lines on
$\Sigma$. Since the general canonical curve is smooth by
assumption, and thus in particular
$|\k|$ has no fixed part, $D$  is a canonical curve.
Notice also that $A^2, B^2\ge 0$ since they move in linear
systems without fixed components, and $AB\ge 3$. On the
other hand,
$AB=A(\k-A)$ is even and therefore $AB\ge 4$. This implies
that $|\k|$ has base points and that there exists a base
point $Q$ of $|\k|$ that is also a base point for both
$|A|$ and
$|B|$. So $A^2, B^2>0$ and $\k^2=A^2+B^2+2AB\ge 10$,
contradicting Claim $1$. We conclude that $\Sigma$ is a
quadric cone.

As in the proof of Claim $2$, consider the
pull-back $D$ of a hyperplane section of $\Sigma$ passing
through  the vertex and write $D=2C+Z$, where $C$ is the pull--back of a
ruling of
the cone  and
$Z$ is an effective divisor, possibly zero, that is contracted to
the vertex. Notice that $C^2\ge 0$, as $C$ moves in a
linear system without fixed components, and $2C^2\le \k C$.
One has
$\k C=DC\ge 3$ and
$2C\k\le
\k^2\le 9$, so that either
$\k C=3$ or $\k C=4$. In the former case, $C^2$ is odd and
thus, by the above discussion,
$C^2=1$. Again (as in the proof of the preceding claim) it
is easy to see that the general curve in $|C|$ is non
hyperelliptic and so $S$ contains a genus 3 pencil as
stated.

 On the other hand, if $\k C=4$  then
$|\k|$ has
 base points and
 there exists a base point $Q$ of $|\k|$ that is also a
base point of
$|C|$, and so $C^2>0$. So we conclude that $C^2=2$ and, by
the index theorem, that
$\k^2\le 8$. But $\k^2\ge 2\k C\ge 8$, so $\k^2=8$. Using
the same argument as in the previous claim it is easy to
see that the general curve in $|C|$  is not hyperelliptic.
\qed

\section{Surfaces with a special pencil $I$}\label{secg3}
\setcounter{defn}{0}
\setcounter{equation}{0}

In this section we study in detail the class of the surfaces containing
a pencil $C$, with $C^2=1$, $\k C=3$, to which one of the
types of surfaces found in section \ref{3-1} belongs. So,
throughout this section, we make the following assumption:

\begin{assu}\label{g3} $S$ is a minimal surface of general
type with $p_g\geq 4$ containing a positive-dimensional
linear system $|C|$ with $\k C=3$, $C^2=1$, and such that
the general curve in $|C|$ is irreducible,
non-hyperelliptic.
\end{assu}

 We describe the possible invariants of such surfaces, show that the canonical
system of these surfaces contains a smooth curve and finally establish
some properties
that will enable us to construct examples of all these surfaces. In
particular we prove in this section  the following:

\begin{thm}\label {I}  Let $S$ be a surface satisfying assumption \ref{g3}.
Then

i) $q=0$ and either $p_g=5$ and $\k^2=9$, or $p_g=4$ and $6\leq \k^2\leq 8$;

ii) $|C|$ is a  pencil with a simple base point $P$;

iii) $\deg\phi_K=3$;

iv) the general curve in $|\k|$ is smooth;

v) the linear system $|3C|$ defines a morphism
$S\to\pp^4$  of degree $3$ onto the normal rational cubic
cone.

\end{thm}

To prove theorem \ref{I} we will need various facts, that
we now establish.

\begin{lem}\label{claim0} If $S$ is a surface as in
\ref{g3}, then $|C|$ is a pencil with a simple base point.
\end{lem}

\proof Remark first of all that $|C|$ has no fixed
components, since it is positive dimensional and the
general curve of $|C|$ is irreducible by assumption. So, if
$\dim|C|$ were greater than $1$, then
$|C|$ would be base point free, and therefore a general
$2$-dimensional subsystem  of
$|C|$ would define an isomorphism of $S$ with $\pp^2$. So
$|C|$ is necessarily
$1$-dimensional and therefore it has a simple  base
point.\qed

 \begin{nota} We  denote by $C$ a general element of $|C|$
and by $P$ the base point of the pencil $|C|$.
\end{nota}
 \begin{lem}\label{claim1} If $S$ is a surface as in assumption
\ref{g3}, then
$\k^2\leq 9$ and
$\k^2=9$ if and only if
$\k\sim 3C$. Furthermore $|\k|$ is not composed with $|C|$.
\end{lem}
\proof By the index theorem, the conditions
 $C^2=1$, $\k C=3$ imply $\k^2\le 9$, with equality
holding if and only if $\k \sim 3C$. By proposition (1.7) of
\cite{cfm}, if $|\k|$ is composed with $|C|$ then
$\k^2=9$, $p_g=4$ and $|\k|=|3C|$. But this  contradicts
\cite{zu}, section 2.6.2.
\qed
\begin{lem}\label{claim2} If $S$ is a surface as in
assumption \ref{g3}, then:

i) $p_g=4$ or $p_g=5$, and $p_g=5$ iff $\k\equiv 3C$;

ii) $h^0(S,\OS(\k-nC))=p_g-n-1$, for $n=1,2,3$;

iii) if $C$ is general, then $|\k|_C$  is a complete and
base point free linear system of dimension $2$.
\end{lem}
 \proof Lemma \ref{claim1} yields $\k(\k-3C)\leq 0$, with
equality holding if and only if $\k\sim 3C$. Thus
$h^0(S,{\cal O}_S(\k-3C))\le 1$, with equality holding if
and only $\k\equiv 3C$.

Consider the restriction maps $$r_n:H^0(S,{\cal
O}_S(\k-nC))\to H^0(C,{\cal O}_C(\k-nC))$$ and notice that
$\ker r_n\simeq H^0(S,{\cal O}_S(\k-(n+1)C))$. Since $\k
C=3$, one has $h^0(C,\OC(\k))\le 2$; on the other hand,
since
$|\k|$ is not composed with $|C|$ by lemma \ref{claim1},
one has $\dim\im r_0\ge 2$.   Thus $h^0(C,\OC(\k)=2$, the
map $r_0$ is onto and
$h^0(S,{\cal O}_S(\k-C))= p_g-2$. Now
$(\k-C)C=2$, hence $\dim\im r_1=h^0(C, \OC(\k-C)=1$,
because $\k-C$ is effective, and $C$ is non--hyperelliptic
by assumption. Thus
$h^0(S,{\cal O}_S(\k-2C))=p_g-3$. Similarly, one gets
$h^0(S,{\cal O}_S(\k-3C))=p_g-4$.
 From $h^0(S,{\cal O}_S(\k-3C))\leq 1$ we have then that
$p_g\le 5$ and
$p_g=5$ iff
$\k-3C\equiv 0$. From the above analysis it follows in
particular that  $|\k|_C$ is a complete linear system of
dimension $2$ and degree $3$. Since $C$ is not
hyperelliptic,
$|\k|_C$ has no fixed point.
\qed
\begin{lem}\label{claim5} If $S$ is a surface as in
assumption \ref{g3} and
$p_g=4$, then
$6\leq \k^2\leq 8$, and $\k\equiv 2C+Z$, where  $Z$ is a
$2$-connected curve with
$ZC=1$, $Z\k=\k^2-6$  and
$Z^2=\k^2-8$. Moreover, if $\k^2=6$  then $Z$ is a smooth
rational curve; if
$\k^2=7,8$  then
$Z$ contains the base point
$P$ of
$|C|$  and there exists
$C'\in |C|$ such that
$Z\subset C'$.
\end{lem}

 \proof Suppose $p_g=4$. By lemma \ref{claim2} ii),
$h^0(S,{\cal O}_S(\k-2C))=1$, and hence
$\k\equiv 2C+Z$, where $Z$ is an effective divisor.  One
has $ZC=1$, since
$\k C=3$, and $\k^2= 2\k C + \k Z\geq 6$, since $\k$ is
nef. It follows immediately that
$Z^2=\k^2-8$ and $\k Z=\k^2-6$. In particular, remark that
if $\k^2=9$ then
$Z\sim C$ by lemma \ref{claim1}.

Let us notice that $Z$ is $2$-connected. In fact assume
otherwise. Then $Z$ decomposes as $Z=A+B$ with $AB\leq 1$.
Since canonical divisors are
$2$--connected, one has
$2\le A(\k-A)=A B+2AC$ and thus
$AC\geq 1$. Similarly we have $B C\ge1$, but this
contradicts
$ZC=1$. Remark that if $Z^2=-2$, $p_a(Z)=0$ and therefore
$2$-connectedness of
$Z$ implies that  $Z$ is an irreducible curve.
 Now let us see that if
$\k^2\geq 7$ then the base point $P$ of $|C|$ lies on $Z$.
In fact, in this case
$Z$ is a $2$-connected curve with $p_a(Z)\geq 1$ and hence
$h^0(Z,{\cal O}_Z(C))=1$ (see \cite{cfm}, (A.5)). This
implies that $|C|$ has a base point lying on $Z$ and thus
$P\in Z$.  Since $P$ lies on $Z$, the restriction map
$$r:H^0(S,{\cal O}_S(C))\to H^0(Z,{\cal O}_Z(C))$$ has
$1$--dimensional image, and therefore there is a curve
$C'\in |C|$ such that $C'=Z+\Delta$, with
$\Delta\geq 0$. Now if $\k^2=9$, then $Z\sim C$ implies
that $\Delta=0$ and
$Z\equiv C$. Therefore,  this case  does not occur by
lemma \ref{claim2}, i).\qed

\begin{lem}\label{claim3}
 If $S$ is a surface as in assumption \ref{g3}, then the
canonical map
$\phi_K$ of
$S$ maps $S$ $3$--$1$ onto the rational normal cone in
$\pp^{p_g-1}$. The curves of $|C|$ are mapped $3$--$1$
onto the rulings of the cone and the base point $P$ is
mapped to the vertex. The map $\phi_K$ is a morphism iff
$\k^2=6$ or
$\k^2=9$.
\end{lem}
\proof Remark that by lemma \ref{claim2}, iii), $\phi_K$
maps a general $C$
$3$--$1$ onto a line. Moreover
$\phi_K$ separates the curves of $|C|$, and thus
$\deg\phi_K=3$.

Let $\Sigma$ be the canonical image of $S$ and let
$d=\deg\Sigma$;  one has:
\begin{equation}\label{grado}
\k^2\ge \deg\phi_K d= 3d\ge 3(p_g-2)
\end{equation}
 By lemmas \ref{claim2}, \ref{claim5} either one has
$\k^2=9$, $p_g=5$, or $6\le \k^2\le 8$, $p_g=4$. In the
former case, (\ref{grado}) implies $d=3$ and $|\k|$ base
point free; in the latter case, it implies
$d=2$ and
$|\k|$ is base point free iff $\k^2=6$. So $\Sigma$ is a
surface of minimal degree, ruled by the images of the
curves of $|C|$.  To show that $\Sigma$ is actually a cone
it is enough to remark that, by lemma \ref{claim2}, iii),
$P$ is not a base point of
$\phi_K$ and so the image lines of the curves of $|C|$ all
go through the point
$\phi_K(P)$.\qed
\begin{lem}\label{claim4} If $S$ is a surface as in
assumption \ref{g3}, then $q=0$.
\end{lem}
  \proof Let us remark first that $q \leq 1$. Indeed
notice that, since
$C^2=1$,   $h^0(C,{\cal O}_{C}(C))=1$  and the restriction
map $H^0(S,{\cal O}_S(C))\to H^0(C,{\cal O}_{C}(C))$ is
onto. Furthermore,
 by lemma
\ref{claim2}, ii), $h^0(C,{\cal O}_{C}(\k))=2$ and  the
restriction map $H^0(S,{\cal O}_S(\k))\to H^0(C,{\cal
O}_{C}(\k))$ is also surjective. So we can choose
$t_0, t_1
\in H^0(S,{\cal O}_S(\k))$ and $s \in H^0(S,{\cal
O}_S(C))$, such that the images of
$t_0, t_1$ generate
$H^0(C,{\cal O}_{C}(\k))$ and the image of $s$ generates
$H^0(C,{\cal O}_{C}(C))$. Then $st_0, st_1$  map to
linearly independent sections of
$H^0(C,{\cal O}_C(\k+C))=H^0(C,\omega_C)$. Consider now the
exact sequence:
 $$0 \to {\cal O}_S(\k)
\to {\cal O}_S(\k+C)
\to \omega_C\to 0.$$ Since $C^2=1$ and $C$ is irreducible,
we have
$h^1(S,{\cal O}_S(\k+C))=0$ (see \cite{bombieri}, pg. 178),
and the long cohomology sequence yields:
$$ H^0(S,{\cal O}_S(\k+C)) \to H^0(C,\omega_{C})\to
H^1(S,\OS(\k)\to 0.$$
 Hence $q=h^1(S,\OS(\k))\le 1$.

 If $q=1$, then  the inequality $\k^2\geq 2p_g$ (see \cite
{debarre}, Th. 6.1) and lemma \ref{claim5} leave us with
the case $p_g=4$, $\k^2=8$. Then
$\chi(\OS)=4$ and
$\k^2<{8\over 3} \chi ({\cal O}_S)$.  By \cite{horikawa-V}, then  the
 Albanese pencil is a genus $2$ pencil. But this
contradicts lemma \ref{claim3}, because the canonical map
of a surface  with a genus $2$ pencil has even degree.\qed
\begin{lem}\label{claim6} If $S$ is a surface as in assumption
\ref{g3}, then $h^1(S,{\cal O}_S(3C))=0$, $h^0(S,{\cal
O}_S(3C))=5$ and $|3C|$ is base point free.
\end{lem}
\proof Notice first that the image of the restriction map
$H^0(S,{\cal O}_S(3C))\to H^0(C,{\cal O}_C(3C))$ has
dimension at most $2$ and the image of
$H^0(S,{\cal O}_S(2C))\to H^0(C,{\cal O}_C(2C))$ has
dimension $1$, since $C$ is not hyperellliptic; so,
arguing as in the proof of lemma \ref{claim2}, it is easy
to show that
$h^0(S,{\cal O}_S(3C))\leq 5$. If $p_g=5$, then by lemma
\ref{claim2}, i),
$h^0(S,{\cal O}_S(3C))=p_g=5$. If $p_g=4$, then
$h^2(S,\OS(3C))=h^0(S,{\cal O}_S(\k-3C))=0$ by lemma
\ref{claim2}, ii). Thus the Riemann-Roch theorem yields
$h^0(S,{\cal O}_S(3C))= \chi ({\cal O}_S)+h^1(S,{\cal
O}_S(3C))=5+h^1(S,{\cal O}_S(3C))$. (Recall that $q=0$ by
lemma \ref{claim4}). So we conclude that
$h^0(S,{\cal O}_S(3C))=5$, $h^1(S,{\cal O}_S(3C))=0$ and
$|3C|_C$ is equal to the complete  linear system $|3P|$.

Finally, observe that the only possible base point of
$|3C|$ is the base point $P$ of
$|C|$; on the other hand if $P$ were a base point of
$|3C|$, then the moving part of the system $|3C|_C$ would
be a $g^1_2$ and $C$ would be hyperelliptic.\qed
\begin{lem}\label{claim7} If $S$ is a surface as in
assumption \ref{g3}, then the linear system $|3C|$ defines
a morphism $f:S\to\pp^4$ of degree $3$ onto the rational
normal cubic cone $C_3$. The curves of $|C|$ are the
pull-backs via $f$ of the rulings of
$C_3$ and the point $P$ is mapped to the vertex.
\end{lem}
\proof Since $C^2=1$ and $|3C|$ is base point free by
lemma \ref{claim6}, every curve of
$|C|$ is mapped by $f$ $3$--to--$1$ onto a line through
the point $f(P)$. So the image of
$S$ is the rational cubic cone and $f$ has degree $3$.\qed

\begin{rem}\label{remarkg3} If $p_g=5$, then the morphisms
$f$ and $\phi_K$ coincide by lemma \ref{claim2}, i). If
$p_g=4$ and $P\in Z$ (see lemma \ref{claim5}), then
$|\k|\subset |3C|$ and
$\phi_K$ is  the composition of $f$ with projection from a
point of $C_3$. Finally, if $\k^2=6$, $p_g=4$ and $P\notin
Z$, then the systems $|3C|$ and
$|\k|$ restrict to different linear systems on a general
$C$ and therefore $f$ and
$\phi_K$ are not related by a birational transformation of
the cones preserving the rulings.
\end{rem}

To finish the proof of thm. \ref{I},  we show that for these surfaces the
general curve in
$|\k|$ is smooth and therefore they satisfy the hypothesis
of theorem \ref{III}.

\begin{lem}\label{claimN} If $S$ is a surface as in assumption \ref{g3}, then
the general curve in $|\k|$ is smooth.
\end{lem}
\proof  For $\k^2=9, 6$, this is immediate since  $|\k|$ is base point free by
lemma \ref{claim3}. For the other cases ($\k^2=7,8$) the canonical map is
$3$--$1$ onto the quadric cone. If $|\k|$ has no fixed part, then its base
locus is a zero--dimensional scheme of length
$\le2$, and therefore it is smooth: thus the general canonical curve is smooth
by Bertini's theorem.  Therefore assume  that $|\k|$ has a fixed part $F$ and
write $\k=M+F$. One has: $M^2\ge 6$, since the linear system $|M|$ maps $S$
$3$--$1$ onto the quadric cone, and
$MF\ge 2$, since  canonical divisors are $2$--connected.  Moreover the
canonical divisor of a minimal surface is nef, and thus we have:
$8\ge\k^2=\k M+\k F\ge \k M=M^2+MF\ge 8$, and thus the
 only  possibility  is $\k^2=8$, $\k F=0$,
$M^2=6$ and $F^2=-2$. Therefore $|M|$ is base point free, and every irreducible
component $\theta$ of $F$ is a smooth rational curve with self intersection
$-2$, i.e a $(-2)$-curve. For any such $\theta$ we have $\theta C=0$, since $C$
is nef and $CF=0$, and $\theta Z\geq 0$ since
$Z$ is $2$--connected by lemma \ref{claim5}.

Consider the effective divisor $\Gamma$ such that $\Gamma$ is linearly
equivalent to $C-Z$ (see lemma \ref{claim5}). One has $\Gamma Z=1$ and $\Gamma
C=0$, hence
$\Gamma^2=-1$, $\k \Gamma=1$. Furthermore $\Gamma$ is 2-connected. In fact,
assume otherwise. Then we can write
$\Gamma=A+B$ where $A,B$ are effective non-zero divisors such that $AB\leq 1$,
$AC=BC=0$. Since $C=\Gamma+Z$ and $\k=2C+Z$, we can write $\k=A+B+C+2Z$: if,
for instance, $AZ=0$ then the divisors $A$ and
$C+2Z+B$ give a decomposition of $\k$ contradicting the fact that $\k$ is
$2$--connected. So we have $AZ>0$ and, by the same argument, $BZ>0$,
contradicting $\Gamma Z=1$.
 Therefore $\Gamma$ is $2$--connected, hence for any irreducible component
$\theta$ of $F$, $\theta \Gamma\geq 0$. Since $\theta C=\theta (Z+\Gamma)=0$,
we have then necessarily $\theta \Gamma=\theta Z=0$ and therefore
$F\Gamma=0$.  But then $M \Gamma=1$,  and so  by proposition (A.5) of
\cite{cfm}), the restriction map
$$r:H^0(S,{\cal O}_S(M))\to H^0(\Gamma,{\cal O}_{\Gamma}(M))$$  has
$1$-dimensional image, a contradiction since $|M|$ is base point free.

Therefore $|\k|$ has no fixed components and thus we proved the lemma.
\qed

\noindent{\bf Proof of Theorem \ref{I}:} The theorem
follows directly from lemmas \ref{claim0},...,
\ref{claimN}.\qed

We continue this section by establishing some facts that will
enable us to construct
 examples of surfaces satisfying assumption \ref{g3} and
having all the possible values of the invariants.
\begin{lem} \label{claim8}  Let $S$ be a surface as in
assumption
\ref{g3}; then

i) if $\k^2\ge7$ or $\k^2=6$ and $P\in Z$ (see lemma
\ref{claim5}), then $h^0(S,
\OS(nC))=6+n(n-3)/2$   for
$n\ge 4$, $P$ is a simple base point of
$|5C|$, and $|nC|$ is  base point free for
$n=4$ and $n\ge 6$;

iii) if $\k^2=6$ and $P\notin Z$ (see lemma \ref{claim5}),
then
 $h^0(S,
\OS(nC))=5+n(n-3)/2$ for $n\ge 4$, $P$ is a simple base
point of
$|4C|$, and $|nC|$ is base point free for
$n\ge 5$.
\end{lem}
\proof

\noindent Case i):

Assume first that $p_g(S)=5$. By lemma \ref{claim2}, for
$n\ge 4$  $nC$ is the adjoint of a nef and big divisor, so
by Kawamata--Viehweg vanishing
$h^0(S,\OS(nC))=\chi(\OS(nC))$ can be computed by
Riemann--Roch.

Assume now that $p_g(S)=4$. By lemma \ref{claim5},
$K_S(K_S-nC)=K_S^2-3n\le 8-3n<0$, and thus
$0=h^0(S,\OS(K_S-nC))=h^2(S,\OS(nC))$  for $n\ge 3$. So we
have:
\begin{equation}\label{RR} h^0(S,
\OS(nC))=5+n(n-3)/2+h^1(S,\OS(nC)),\quad n\ge 3.
\end{equation} Denote by $C$ a smooth element of $|C|$ and
consider the sequence:
\begin{equation}\label{sequence}
0\to\OS\to\OS(C)\to\OO_C(P)\to 0
\end{equation}
 Twisting (\ref{sequence}) by $3C$ and recalling that
$h^1(S,\OS(3C))=0$ by lemma
\ref{claim6}, one obtains  exact sequences:
$$0\to H^0(S,{\cal O}_S(3C))\to H^0(S,{\cal O}_S(4C))\to
H^0(C,{\cal O}_C(4P))\to 0$$ and
$0\to H^1(S,\OS(4C))\to H^1(C,\OC(4C))\to 0$. By the
adjunction formula, $4P$ is a canonical divisor on $C$ and
thus
 $h^1(S,\OS(4C))=1$, $h^0(S,\OS(4C))=7$. Twisting
(\ref{sequence}) by $nC$, $n\ge 4$, and passing to the
associate cohomology sequence one obtains a surjection
$H^1(S,\OS(nC))\to H^1(S,\OS((n+1)C))$, since
$h^1(C,\OO_C((n+1)P))=0$. So we have  $h^1(S,\OS(nC))\le
1$ for $n\ge 4$. By Serre duality, $h^1(S,\OS(nC))=h^1(S,
\OS(-(n-3)C-\Gamma))$, where
$\Gamma$ is the only effective divisor linearly equivalent
to $C-Z$ (cf. lemma
\ref{claim5}). Since $\Gamma C=0$, if $D$ is a general
element of
$|(n-3)C|$ then
$D+\Gamma$ is a disconnected curve and
$h^0(D+\Gamma,\OO_{D+\Gamma})\ge 2$. This implies that
$h^1(S,\OS(nC))=h^1(S,
\OS(-(n-3)C-\Gamma))\ge 1$ and so, eventually,
$h^1(S,\OS(nC))=1$ for $n\ge 4$. So $h^0(S,\OS(nC))$ can
now be computed from (\ref{RR}). Using again Riemann-Roch
on the curve $C$, one sees that $|5P|=P+|4P|$ and thus $P$
is a base point of $|5C|$. It is a simple base point (and
the only one), since
$|4C|$ is free and $P$ is the only base point of $|C|$. We
have shown above that for $n\ge 4$ the map
$H^1(S,\OS(nC))\to H^1(S,\OS((n+1)C))$ is an isomorphism;
it follows that $H^0(S, \OS(nC))\to H^0(C,\OC(nP))$ is
surjective for $n\ge 5$ and thus $|nC|$ is free for $n\ge
6$.
\smallskip

 \noindent Case ii):

The arguments here are analogous to
those used in the proof of case i). Notice that
$\k(\k-nC)=6-3n<0$ for $n\ge 3$, so that
$h^2(S,\OS(nC))=h^0(S,
\OS(\k-nC)=0$ for $n\ge 3$.  Fix a smooth
$C\in |C|$ and denote by
$Q$ the intersection point of
$C$ and
$Z$. We have $P\ne Q$ by assumption. By the adjunction
formula, $3P+Q$ is a canonical divisor on $C$ and so,
using Riemann--Roch on $C$, one shows easily that
$|4P|=P+|3P|$. Sequence (\ref{sequence}) and the fact
that, by lemma
\ref{claim6},
$h^1(S,\OS(3C))=0$  imply that $P$ is a base point of
$|4C|$, $h^0(S,\OS(4C))=7$ and $h^1(S,\OS(4C))=0$. Since
$|3C|$ is free and
$P$ is a simple base point of $|C|$, $P$ is also a simple
base point of $|4C|$. Furthermore,   for $n\ge 4$ the map
$H^1(S,\OS(nC))\to H^1(S,\OS((n+1)C))$ is onto, since
$H^1(C,\OO_C((n+1)P))=0$. So $h^1(S,\OS(nC))=0$ for $n\ge
4$, $H^0(S, \OS(nC))\to H^0(C, \OC(nP))$ is onto for $n\ge
5$ and thus $|nC|$ is free for $n\ge 4$. Finally,
$h^0(S,\OS(nC))$ can now be computed by means of the
Riemann--Roch formula.\qed

\section{Surfaces with a special pencil $I\!I$}
\setcounter{defn}{0}
\setcounter{equation}{0}

In this section we study the surfaces containing a
positive-dimensional linear system $|C|$ satisfying
$K_SC=4$, $C^2=2$ and such that  the general curve in
$|C|$ is irreducible, non-hyperelliptic, to which belong
one of the types of surfaces we encountered in section
\ref{3-1}.

  So, throughout this section, we make the following
assumption
\begin{assu}\label{g4}  $S$ is a minimal surface of
general type with
$p_g\ge 4$, containing a positive--dimensional linear
system $|C|$ satisfying $K_SC=4$, $C^2=2$, and such that
the general curve in $|C|$ is irreducible,
non--hyperelliptic.
\end{assu}
\begin{nota} We will denote by $C$ a general element of
$|C|$.
\end{nota}

We want to prove the following theorem and also to
establish some properties that will enable us to construct
examples.

\begin{thm}\label {II} Let $S$ be a surface as in
assumption
\ref{g4}. Then

i)  $\k^2=8$, $p_g=4$ and $q=0$;

ii) $|C|$ is a  pencil with two base points $P$ and
$P_2$;  $\k=2C$ and the general  curve of $\k$ is smooth;

iii) $\deg\phi_K=3$ or
$\deg\phi_K=4$;

iv)  if
$\deg\phi_K=3$,  $P_2$ is infinitely near to $P$; if we denote
by $S'$ the blow-up of $S$ at $P$  and by $C'$ the strict
transform of $C$ on $S'$, then the linear system
$|3C'|$ defines a morphism $S\to\pp^4$  of degree $3$ onto
the normal rational cubic cone.
\end{thm}
To prove the theorem we will need the following:
\begin{lem}\label{1} If $S$ is a surface as in assumption
\ref{g4}, then
$|C|$ is a pencil,
$p_g=4$, $\k^2=8$, $\k=2C$, $\phi_K$ is generically finite,
 $\deg\phi_K=3$ or
$4$ and the general curve of $|\k|$ is smooth .
\end{lem}
\proof
 Notice first that the assumptions imply that the general
curve in $|C|$ is smooth irreducible of genus $4$.
Moreover, if $\dim|C|>1$, then the moving part of the
restriction of $|C|$ to a general $|C|$ would be a $g^1_2$
and $C$ would be hyperelliptic. So $|C|$ is a pencil. For
a general
$C\in|C|$, one has, by Riemann--Roch:
$h^0(C,\OC(\k))=h^0(C,\OC(C))+1=2$, the last equality
holding because $C$ is not hyperelliptic. Now, arguing as
in the proof of lemma \ref{claim2}, one shows that
$h^0(S, \k-2C)\ge p_g-3\ge 1$. So, $\k^2-8=\k(\k-2C)\ge
0$,  since $\k$ is nef. On the other hand,  the index
theorem gives $\k^2\le 8$, with equality holding if and
only if $\k\sim 2C$. Since we have shown that $\k-2C$ is
effective, we conclude that $\k=2C$ and $\k^2=8$. In particular,
$1=h^0(S, \OS(\k-2C))\ge p_g-3\ge 1$ and so
$p_g=4$.

The image of $\phi_K$ is a surface, since otherwise
$\phi_K$ would be composed with $C$ and thus $p_g(S)=3$.
The assertion about the degree of $\phi_K$ can be proven
exactly as in lemma
\ref{claim3} and the last assertion is obvious from $\k=2C$. \qed
\begin{lem} If $S$ is a surface as in assertion \ref{g4},
then $q=0$.
\end{lem}
\proof  This lemma is proved for a slightly different
situation in
\cite{cfm}, proposition (3.5). For the reader's
convenience we give here an outline of the proof.

First of all, arguing as in the proof of lemma
\ref{claim4} one shows  that
$q\le 2$. Next assume that $q\neq 0$ and let $\mu$ be an
element in $Pic^0(S)$. From the exact sequence
$$ 0 \to \mu \to \mu(C) \to {\cal O}_{C}(C)\otimes \mu \to
0$$ we see that, for $\mu$ general, one has
$h^0(S,\mu(C))\leq 1$,
 because
$C^2=2$.  By the same reason, also
$h^0(S,\mu^{\vee}(C))\leq 1$. By Serre duality and lemma \ref{1}, it follows
$h^2(S,\mu(C))=h^0(S,\mu^{\vee}(C)$ and thus, by the Riemann-Roch theorem, one
has:
\begin{equation} 2\ge h^0(S,\mu(C)) +
h^0(S,\mu^{\vee}(C))=p_g-q+ h^1(S,\mu(C))\ge p_g-q
\end{equation}
Since $p_g=4$ and $q\le 2$, the above
inequality yields $q=2$ and $h^0(S,
\mu(C)) = h^0(S,\mu^{\vee}(C))=1$.

So we can define a rational map $\rho: Pic^0(S) \to
C^{(2)}$, which sends  a general element $\mu$ to the only
effective divisor in
$|\OC(C)\otimes
\mu|$. The map
$\rho$ is generically injective. Since $C$ is a smooth
curve of genus $4$, we find a contradiction, because the
surface of general type $C^{(2)}$ cannot be dominated by
the abelian surface
$Pic^0(S)$.
\qed
\begin{lem}\label{basepointg4} Let $S$ be a surface as in
assumption \ref{g4}, and assume that $\deg
\phi_K=3$; then:

i) the image of $\phi_K$ is a quadric cone;

ii) the  base locus of $|\k|$ consists of a point $P$  and
a point
$P_1$, infinitely near to $P$ and in particular the general canonical curve is
smooth;

iii) the base locus of
$|C|$ consists of $P$ and a point $P_2$ infinitely near to
$P$, with $P_2\ne P_1$.
\end{lem}
\proof Since $\k\equiv 2C$, $|\k|$ has no fixed component.
Denote by $d$ the degree of the image $\Sigma$ of
$\phi_K$; one has
$8=\k^2\ge d\deg\phi_K=3d\ge 3(p_g-2)=6$. So $d=2$ is the
only possibility and
$\Sigma$ is a singular quadric, since $\k=2C$. Moreover
the above inequality implies that $|\k|$ has two simple
base points
$P$ and
$P_1$, with
$P_1$  possibly infinitely near to $P$. Notice that by Bertini's theorem this
implies that the general canonical curve is smooth.
 We wish to show
that
$P_1$ is actually infinitely near to $P$. So assume
otherwise: then both $P$ and
$P_1$ are base points of $|C|$ and thus also of
$|\k|_C$.  This implies that the moving part of $|\k|_C$
is a $g^1_2$ and
$C$ is hyperelliptic, against the assumptions. Notice that
$P$ is necessarily a base point of $C$; since $C^2=2$, $C$
has also another base point $P_2$. We are going to prove
that $P_2$ is also  infinitely near to $P$ but $P_2\ne
P_1$. By adjunction, the sheaf
$\OS(2\k)$ restricts to $\omega_C(P+P_2)$ on $C$; by
Riemann--Roch on $C$, this implies that $|2\k|$ does not
separate $P$ and
$P_2$, and thus also $|\k|$ does not separate $P$ and
$P_2$. If $P_2$ is not infinitely near to $P$, this is
equivalent to saying that $P_2$ is also a base point of
$|\k|$, contradicting what we have just proven. So we
conclude that $P_2$ is infinitely near to $P$. Finally, if
$P_1$ and $P_2$ were equal, then again the moving part of
$|\k|_C$ would be a
$g^1_2$, and $C$ would be hyperelliptic.\qed
\begin{nota}\label{blowup} We denote by $\epsi':S'\to S$
the blow--up of $S$ at the common base point $P$ of $|C|$
and $|\k|$,  by $E'$ the exceptional curve of $\epsi'$,
again by
$C$ the pull--back of $C$ to $S'$ and by $C'$ the divisor
$C-E'$.
\end{nota}
\begin{lem}\label{freeg4} Let $S$ be a surface as in
assumption \ref{g4} such that $\deg\phi_K=3$, and let
$S'$, $C'$ be as in notation \ref{blowup}; then
$h^0(S',\OO_{S'}(3C'))=5$ and the system $|3C'|$ is base
point free.
\end{lem}
\proof By lemma \ref{basepointg4}, the moving part of the
canonical system of $S'$ is
$|M|=|2C'+E'|$. Let $C'\in|C'|$ be a general curve: then
$|M|_{C'}=|3P_2|$ is a base point free complete linear
system of dimension $1$ and degree $3$, since the
canonical map has degree $3$. Using restriction sequences
as in the proof of lemma \ref{claim2} it is easy to show
that $h^0(S',\OO_{S'}(2C'))=3$ and that
$h^0(S',\OO_{S'}(C'-E'))=1$. Thus
$|M|\subset |3C'|$ and  $|3C'|$ is not composed with
$|C'|$. From the above considerations it follows that the
sequence:
$$0\to H^0(S',\OO_{S'}(2C'))\to H^0(S',\OO_{S'}(3C'))\to
H^0(C',\OO_{C'}(3P_2))\to 0$$ is exact. This amounts to
saying that $h^0(S',\OO_{S'}(3C'))=5$ and $P_2$ is not a
base point of $|3C'|$. Since $P_2$ is the only base point
of $|C'|$ by lemma
\ref{basepointg4}, it follows that $|3C'|$ is free. \qed
\begin{lem}\label{tripleg4}
 If $S$ is a surface as in assumption \ref{g4} and $S'$,
$C'$ are as in notation
\ref{blowup}, then the linear system $|3C'|$ defines a
morphism $f:S\to\pp^4$ of degree
$3$ onto the normal rational cubic cone
$C_3$. The curves of $|C'|$ are the pull--backs via $f$ of
the rulings of
$C_3$ and the point $P_2$ (see lemma \ref{basepointg4}) is
mapped to the vertex.
\end{lem}
\proof Lemma \ref{freeg4} and its proof show that $f$ maps a
general curve of $|C'|$  $3$--to--$1$ onto a line in $\pp^4$ contaning the
point $f(P_2)$. So the image of $f$ is a cone with vertex
$f(P_2)$. Since $(3C')^2=9$, the image of $f$ has degree
$3$, and so it is the normal rational cubic cone. \qed
\noindent {\bf Proof of Theorem \ref{II}:} The theorem
follows directly from lemmas
\ref{1},\ldots \ref{tripleg4}.\qed
We close this section
by establishing some results on linear systems on
$S'$ that that will be needed for the contruction of
examples of surfaces as in assumption
\ref{g4} with canonical map of degree $3$.
\begin{lem}\label{dimg4} Let $S$ be a surface as in
assumption \ref{g4} such that $\deg\phi_K=3$, and let
$S'$, $C'$, $E'$ be as in notation \ref{blowup}; then:

i) $h^0(S',\OO_{S'}(4C'))=7$, the point $P_2$ (see lemma
\ref{basepointg4}) is the only base point of
$|4C'|$ and it is a simple one;

ii) $h^0(S',\OO_{S'}(5C'))=10$ and $|5C'|$ is base point
free;

iii) if $n\ge 6$, then
$h^0(S',\OO_{S'}(nC'))=11+\frac{n(n-5)}{2}$ and
$|nC'|$ is base point free for $n\ne 7$.
\end{lem}
\proof Restricting $\OO_{S'}(4C')$ to a general $C'$ and
taking global sections one obtains the sequence:
\begin{equation}\label{sopra} 0\to
H^0(S',\OO_{S'}(3C'))\to H^0(S',\OO_{S'}(4C'))\to
H^0(C',\OO_{C'}(4P_2))\to 0.
\end{equation} By lemma \ref{freeg4}, the system $|3P_2|$
on $C$ is free of dimension $2$ and it is equal to the
restriction of $|3C'|$ to
$C'$. Using Riemann-Roch on $C'$ and the fact that $C'$ is
not hyperelliptic, one shows easily that
$|4P_2|=P_2+|3P_2|$. So sequence \ref{sopra} is right
exact and claim i) follows now easily from lemma
\ref{freeg4}.

As we have already remarked in the proof of lemma
\ref{freeg4}, there exists a unique effective divisor
$C_0\equiv C'-E'$.  By lemma  \ref{basepointg4}, the
canonical map $\phi'_K$ of
$S'$ is induced by the system $|M|=|2C'+E'|$ that has
$P_1$ as its only base point. Since $MC_0=E'C_0=2$ and
$p_a(C_0)=3$, it follows that $P_1\notin C_0$ and
$\phi'_K$ maps
$C_0$
$2$--to--$1$ onto a ruling of the quadric cone in $\pp^3$.
One has $C'C_0=0$ and so the restriction of $\OO_{S'}(C')$
to $C_0$ is trivial. The sequence
$0\to\OO_{S'}(4C'+E')\to\OO_{S'}(5C')\to
\OO_{C_0}\to 0$ induces on cohomology the exact sequence:
$$0\to H^0(S',\OO_{S'}(4C'+E))\to H^0(S',\OO_{S'}(5C'))\to
H^0(C_0,
\OO_{C_0}).$$ Since the image of $H^0(S',\OO_{S'}(5C'))\to
H^0(C_0,
\OO_{C_0})$ is nonzero,  we have
$h^0(S',\OO_{S'}(5C'))\ge h^0(S',\OO_{S'}(4C'+E))+1$; to
compute
$h^0(S',\OO_{S'}(4C'+E'))$ we restrict again to $C_0$ and
consider the sequence on global sections:
$$0\to H^0(S',\OO_{S'}(3C'+2E'))\to
H^0(S',\OO_{S'}(4C'+E')) \to H^0(C_0,\OO_{C_0}(E'))\to 0.$$
By the above considerations,  $|\OO_{C_0}(E')|$ is the
restriction of $|M|\subset |4C'+E'|$ to $C_0$ and is a
$g^1_2$. So the sequence is  exact and
$h^0(S',\OO_{S'}(4C'+E'))=h^0(S',\OO_{S'}(3C'+2E'))+2$.

The last step is now the computation of
$h^0(S',\OO_{S'}(3C'+2E'))$: the divisor
$3C=\k+C$ on the surface
$S$  is the adjoint of a nef and big divisor, so by
Kawamata--Viehweg vanishing $h^0(S,
\OS(3C))=\chi(\OS(3C))=8$. By the regularity of
$S$, $|3C|$ restricts to the complete system $|K_C|$ on a
generic $C$ and thus
$P$ is not a base point of $|3C|$. It follows
$h^0(S',\OO_{S'}(3C'+2E))=7$, and, finally,
$h^0(S',\OO_{S'}(5C'))\ge 10$. Now, taking into account
the dimensions of the vector spaces involved, one sees
that the sequence:
$$0\to H^0(S',\OO_{S'}(4C'))\to H^0(S',\OO_{S'}(5C'))\to
H^0(C',\OO_{C'}(5P_2))\to 0$$ is exact and therefore
$|5C'|$ is base point free of dimension $9$.

If $n\ge 5$, then
$nC'=K_{S'}+(n-2)C'-3E=K_{S'}+(n-5)C'+3C_0$ is the adjoint
of an effective divisor, and thus
$h^2(S',\OO_{S'}(nC'))=0$ and
$\chi(\OO_{S'}(nC'))=h^0(S',\OO_{S'}(nC'))-h^1(S',\OO_{S'}(nC'))=5+\frac{n(n
-5)}{2}$. By the regularity of $S'$, the system
$|K_{S'}+C'|$ restricts to the complete canonical system
$|6P_2|$  on every curve $C'$. One has $6C'=
K_{S'}+C'+3C_0$ and so the restriction map $H^0(S',
\OO_{S'}(6C'))\to H^0(C', \OO_{C'}(6P_2))$ is surjective,
implying $h^0(S', \OO_{S'}(6C'))=h^0(S',
\OO_{S'}(5C'))+h^0(C',
\OO_{C'}(6P_2))=14$ and, as a consequence, $h^1(S',
\OO_{S'}(6C'))=6$. Since $h^1(C',\OO_{C'}(nP_2))=0$ for
$n\ge 7$, one has a surjection $H^1(S',
\OO_{S'}((n-1)C'))\to H^1(S', \OO_{S'}((nC'))$ for $n\ge
7$. So for $n\ge 6$ one has $h^1(S', \OO_{S'}((nC'))\le 6$
and $h^0(S', \OO_{S'}((nC'))\ge 11+\frac{n(n-5)}{2}$, with
equality holding for $n=6$. In order to show that equality
actually holds for $n\ge 7$, it is enough to prove that
the restriction map $H^0(S', \OO_{S'}(nC'))\to H^0(C',
\OO_{C'}(nP_2))$ is surjective for
$n\ge 7$. Using the above discussion, the Riemann--Roch
theorem on $C'$ and the base point free pencil trick, it
is not difficult to prove that the graded ring
$\oplus_{n\ge 0}H^0(C',\OO_{C'}(nP_2))$ is generated by
elements of degree $\le 5$. Since we already know that the
restriction map $H^0(S', \OO_{S'}(nC'))\to H^0(C',
\OO_{C'}(nP_2))$ is surjective for $n\le 5$, it follows
that it is indeed surjective for all $n$.\qed

\section{Triple covers of $\f$}
\setcounter{defn}{0}
\setcounter{equation}{0}

In the previous sections we have shown that all surfaces
satisfying assumption
\ref{g3} and all surfaces satisfying  assumption \ref{g4}
and having canonical map of degree $3$ are generically
finite triple covers of the rational cubic cone
$C_3$ in
$\pp^4$. In this section we describe these triple covers
more precisely, in order to be able to show the existence
of the surfaces we are studying.

Triple covers have been studied in \cite{mi} and
\cite{pa} under the assumption that the map be flat. So we
wish to reduce ourselves to the case of a flat map. We
state first some general results on coverings that will be
useful later:
\begin{lem}\label{flatness} Let $X$, $Y$ be  surfaces,
with $X$ normal and $Y$ smooth and let $f:X\to Y$ be a
quasi--finite projective map: then $f$ is flat.
\end{lem}
\proof The surface $X$, being normal, is Cohen--Macaulay
by Theorem 8.22 A of
\cite{ha}; moreover the map $f$ is finite by Exercise
11.2, chapter $I\!I\!I$ of
\cite{ha}. So we may apply Corollary 18.17 of
\cite{eisenbud} to conclude that $X$ is flat over $Y$.
\qed
\begin{lem}\label{singtriplo} Let $Y$ be a smooth variety
and let $\pi:X\to Y$ be a finite flat map of degree $3$
with branch locus $D$, and let $Q$ be a (not necessarily
closed) point of $Y$:

i) if $\pi$ is simply ramified over $Q$, then
$Q\in D$ and $X$ is smooth over $Q$ iff $D$ is smooth at
$Q$;

ii) if $\pi$ is totally ramified over $Q$, then $D$ is
singular at
$Q$ and
$X$ is smooth over $Q$ iff  $D$ has a double point at $Q$.
\end{lem}
\proof This is just a more synthetic formulation of the
analysis of section
$5$ of
\cite{mi}.
\qed
\begin{lem} \label{normalita'} Let $X$, $Y$ be irreducible
varieties, with $Y$ normal, and let $f:X\to Y$ be a finite
flat map. Then $X$ is normal iff it is nonsingular in
codimension $1$.
\end{lem} {\bf Proof:} It is well known (see \cite{ha},
Thm. 8.22A) that if
$X$ is normal, then it is nonsingular in codimension $1$.

 Assume now that $X$ is not normal. Then there exists an
affine open subset
$U=\spec B$ of $X$ that is not normal. Since
$f$ is finite (and thus affine), we may assume that
$U=f^{-1}V$, where
$V=\spec A$ is an affine open subset of $Y$. By the
flatness of $f$, possibly after shrinking
$V$, there exist
$z_1\ldots z_k\in B$ such that $B=A\oplus Az_1\oplus\cdots
Az_k$. Notice that, if
$F$ is the fraction field of $A$, then the fraction field
of $B$ is
$G=B\otimes_AF$.
 Let $x\in G$ be integral over $B$, but not an element of
$B$; then
$x=a_0+a_1z_1+\cdots a_kz_k$ with at least one of the
$a_i$'s, say $a_j$, not in
$A$. Then the locus where
$a_j$ is not regular is a divisor in $V$, $Y$ being
normal, and thus the set of non--normal points of $X$ has
codimension
$1$. But this contradicts the assumption that $X$ be
nonsingular in codimension
$1$.\qed
\begin{cor}\label{triplenormal} Let $Y$ be a smooth
surface and $\pi:X\to Y$ a finite flat  map of degree
$3$; assume that $X$ is not normal and denote by
$\tilde{X}$ its normalization: then the map
$\tilde{X}\to Y$ is flat and finite, with branch locus
strictly contained in the branch locus  of $X\to Y$.
\end{cor}
\proof The map $\tilde{X}\to Y$ is flat and finite, by
lemma
\ref{flatness}.  Denote by
$D$ the branch locus of
$X\to Y$ and by
$\tilde{D}$ the branch locus of
$\tilde{X}\to Y$. By the flatness of $\pi$ and $f$, $D$ and
$\tilde{D}$ are divisors and
$\tilde{D}\le D$. By  lemma \ref{normalita'}, if
$X$ is not normal then it is singular in codimension 1 and
thus, by lemma
\ref{singtriplo}, there is an irreducible component $C$ of
$D$ such that either $X$ is not totally ramified over $C$
and $C$ appears in $D$ with multiplicity
$\ge 2$ or $X$ is totally ramified over $C$ and $C$
appears in $D$ with multiplicity $\ge 3$. In either case
$C$ appears in $\tilde{D}$ with multiplicity strictly
smaller than in  $D$.
\qed

\begin{nota} We denote by $\f$ the ruled surface
$\proj(\OO_{\pp^1}\oplus\OO_{\pp^1}(3))$, by $\so$ a
tautological section, by $\si$ the infinity section and by
$R$ the class of a ruling of $\f$. So we have
$\so\equiv \si+3R$, $\so^2=3$,
$\si^2=-3$,
$R^2=0$. Recall that the linear system $|\so|$ defines a
morphism
$\psi:\f\to\pp^4$ that maps $\f$ birationally onto the
cubic cone $C_3$ and contracts $\si$ to the vertex $v$ of
$C_3$.
\end{nota}
\begin{prop}\label{morfismo} Let $S$ be as in assumption
\ref{g3}, let $C_3$ be the rational normal cubic cone in
$\pp^4$, let $f:S\to C_3$ be the map induced by the
linear system
$|3C|$ (see lemma
\ref{claim7}), let
$\epsi:\hS\to S$ be the blow-up of
$S$ at the base point
$P$ of $|C|$ (see lemma \ref{claim0}); then there is a commutative diagram:
\begin{equation}
\begin{array}{rcccl}
\phantom{1} &\hS&\stackrel{\epsi}{\rightarrow} &S&
\phantom{1} \\
\scriptstyle{g}\!\!\!\!\!\! & \downarrow & \phantom{1} &
\downarrow &
\!\!\!\!\!\!
\scriptstyle{f}
\\
\phantom{1} & \f & \stackrel{\psi}{\rightarrow} & C_3 &
\phantom{1}
\end{array}
\end{equation} where $g:\hS\to\f$ is a generically finite
morphism of degree $3$. The map
$g$ is totally ramified and finite over
$\si$.
\end{prop} Before proving the proposition we fix some
notation:
\begin{nota} We denote by $E$ the exceptional curve of
$\epsi:\hS\to S$; we still denote by $|C|$ the pull--back
to $\hS$ of the system $|C|$ on $\hS$, and by $|F|$ the
strict transform of $|C|$, so that $|F|$ defines a
fibration $\hS\to\pp^1$.
\end{nota} {\bf Proof of prop. \ref{morfismo}:} In
principle,  $g$ is a rational map of degree $3$; in order
to show that $g$ is a morphism we are going to prove that
the pull--back on $\hS$ of the very ample linear system
$|R+\so|$ is base point free. By lemma \ref{claim7}, one
has:
$g^*R=F=C-E$ and
$g^*\so=3C$. By lemma \ref{claim8}, the system
$|4C-E|=|3C+F|$ has dimension $6 =\dim|R+\so|$,  and
therefore it is equal to
$g^*|R+\so|$. The system $|F|$ is free, since it is
irreducible and $F^2=0$, and the system
$|3C|$ is free by lemma \ref{claim7}, so that $|3C+F|$ is
also free and $g$ is a morphism. The inverse image of
$\si$ contains the divisor $3E$, since $3C|_C=3P$: since
$\si\equiv\so-3R$, then $g^*\si\equiv g^*(\so-3R)\equiv
3C-3R \equiv 3E$, one gets immediately $g^*\si=3E$ and thus
$g$ is totally ramified over $\si$. \qed
\begin{nota-rem}\label{stein} We denote by
$\hS\stackrel{h}{\to}X\stackrel{\pi}{\to}\f$  the Stein
factorization of
$g:\hS\to \f$ (see lemma \ref{morfismo}), so that $X$ is a normal variety,
$h$ is
a morphism with connected fibres and
$\pi:X\to\f$ is a finite  morphism of degree $3$. Notice
that by lemma
\ref{flatness}, $\pi$ is flat.  So  the theory of \cite{mi}
and \cite{pa} applies and
$\pi_*\OX=\OO_{\f}\oplus \E$, where $\E$ is a locally free
sheaf of rank $2$, the so-called {\em trace-zero module}
of $\pi$.
\end{nota-rem} We determine $\E$ in the next proposition.
\begin{prop}\label{Eg3} Let $S$ be a surface as in
assumption \ref{g3} and  let $\pi:X\to
\f$ and ${\cal E}$ be as in notation--remark \ref{stein}.
Then

i) if $\k^2=7,8,9$ or if $\k^2=6$ and $P\in Z$ (see lemma
\ref{claim5}), then $\E=
\Of(-2\si-4R)\oplus
\Of(-3\si-8R)$;

ii) if $\k^2=6$ and $P\notin Z$ , then
$\E=\Of(-2\si-5R)\oplus
\Of(-3\si-7R)$.
\end{prop}

\proof Case i): (For $\k^2=9$  this is the same as
\cite{Konno}, Thm. 2.3).

\noindent   Notice that, by lemma \ref{claim7}, $g^*R=F$
and, by prop.
\ref{morfismo},
$g^*\si=3E$; so if one lets
$L=\Of(2\si+4R)$ and
$M=\Of(3\si+8R)$, one has $g^*L= \OSh(4F+6E)=\OSh(4C+2E)$
and
$g^*M=\OSh(8F+9E)=\OSh(8C+E)$. Since $(g^*L) E=-2$, the
divisor $2E$ is contained in the fixed  part of $|g^*L|$;
on the other hand, by lemma
\ref{claim8}, $|4C|$ is base point free of dimension $7$.
So
$|g^*L|=2E+|4C|$, and $g^*|L|\subset |g^*L| $ is the
$6$--dimensional subsystem of curves vanishing on
$E$ of order $\ge3$.
 Fix a generic $z\in H^0(\hS,g^*L)$: then $z$ vanishes on
$E$ of order $2$ and is not a pull-back from $\f$. Arguing
as above one sees that $|g^*M|=E+|8C|$, where $|8C|$ is
free, since $|4C|$ is, while $g^*|M|$ consists of divisors
vanishing on $E$ of order $\ge 3$.  Again we may choose
$w\in
\Ho(\hS,g^*M)$ that vanishes on $E$ precisely of order $1$
and is not a pull--back from $\f$. The pair
$(z,w)$ defines a birational morphism $\psi:\hS\to L\oplus
M$ such that $g$ factors through $\psi$.  If we can show
that the image $Y$
 of
$\psi$ is a normal surface and the map $Y\to \f$ is
finite, then it will follow that
$\psi:\hS\to Y$ coincides with $h:\hS\to X$ by the
universal property of the Stein factorization.
 We start by determining the  equations defining
$Y$ inside the vector bundle $ L\oplus M$.
 We claim that
$H^0(\hS,g^*L^2)$ is the direct sum of the following
subspaces:
$$ V_0=g^*H^0(\f,L^2),\quad V_1=wg^*\Ho(\f,L^2\otimes
M\inv).\quad V_2=zg^*H^0(\f,L).$$ Notice that, since
$g^*L^2=\OSh(4E+8C)$ and $EC=0$,
$|g^*L^2|=4E+|8C|$. So $h^0(\hS,g^*L^2)=26$, by lemma
\ref{claim8}, and it is easy to check that $\dim V_0=18$,
$\dim V_1=1$, $\dim V_2=7$. Thus it is enough to show that
for any relation of the form $s_0+s_1+s_2=0$, with $s_i\in
V_i$, one has
$s_0=s_1=s_2=0$. This follows at once by remarking  that
if $s_i\ne 0$ then it vanishes on $E$ with order $i\mod 3$,
for $i=0,1,2$.
 So  $z^2\in H^0(\hS,g^*L^2)$ can be written uniquely as
$z^2=az+bw+A$, where
$a\in g^*H^0(\f,L)$, $b\in g^*H^0(\f, L^2\otimes M\inv)$
and $A\in g^*H^0(\f, L^2)$. Computations of the same kind
yield relations of the form
$zw=ez+fw-B$ in
$\Ho(\hS,g^*L\otimes g^*M)$ and $w^2=cz+dw+C$ in
$\Ho(\hS,g^*M^2)$, where $e\in g^*H^0(\f, M)$, $f\in
H^0(\f,L)$, $B\in g^*H^0(\f, L\otimes M)$, $c\in g^*H^0(\f,
M^2\otimes L\inv)$, $d\in g^*H^0(\f, M)$ and $C\in
g^*H^0(\f,M^2)$. If we replace
$z$ by
$z-\frac{a+f}{3}$ and $w$ by $w-\frac{d+e}{3}$, the
relations take the form:
\begin{equation}\label{equazioni}
\left\{
\begin{array}{lcc}  z^2 &=&az+bw+A\\ zw&=&-dz-aw-B\\
w^2&=&cz+dw+C
\end{array}\right.
\end{equation}
 We now wish to show that $A=2(a^2-bd)$, $B=ad-bc$ and
$C=2(d^2-ac)$, so that the equations \ref{equazioni} are
of the same form of those given in Thm. 2.7 of
\cite{mi}. To do this, one computes
$z^2w$ both from the first and the second equation and
equates the two expressions thus obtained; this yields a
relation involving $z^2, zw, w^2, z, w$, from which the
terms of degree $2$ in $w$, $z$ can be eliminated using
equations  \ref{equazioni} (See lemma 2.4 and lemma 2.6
of \cite{mi}). So one ends up with a relation of the form
$sz+tw+r=0$ in
$\Ho(\hS, g^*(L^2\otimes M)$, where $s\in
g^*H^0(\f,L\otimes M)$,
$t\in g^*H^0(\f,L^2)$ and $r\in g^*H^0(\f,L^2\otimes M)$.
Again by considering the vanishing orders on
$E$ of the three summands, one deduces $s=t=r=0$. This
yields
$A=2(a^2-bd), B=bc-ad$. The same elimination procedure for
$zw^2$ finally gives
$C=2(d^2-ac)$.

Equations \ref{equazioni} define in $L\oplus M$ a flat
finite triple cover
$Y'\to\f$ with trace zero module $\E=L\inv\oplus M\inv$,
and $Y$ is contained in
$Y'$.
 Let $U\subset \f$ be the set of regular values of $g$:
$Y|_U$ is an  \'etale cover of degree $3$ and so
$Y|_U=Y'|_U$. This shows that $Y'$ is generically
reduced;  then,   by prop. $3.4$ of
\cite{pa}, it is reduced,  and therefore $Y'=Y$. Next we
are going to show that
$Y$ is normal.  By the universal property of the Stein factorization, there
exists a normalization map $\nu:X\to Y$ such that $\pi$ is the
composition of $Y\to\f$ with $\nu$. Denote by $D$ the
branch locus of $Y\to\f$ and by $D'$ the branch locus of
$\pi:X\to\f$: by cor. \ref{triplenormal}
 $D'\le D$ and, moreover, $D=D'$ if and only if $Y$ is
normal (and in that case $\nu:X\to Y$ is an isomorphism). By prop.
$7$ of
\cite{mi}, $D$ is linearly equivalent to $2\si+8\so$. By
prop.
\ref{morfismo} and lemma \ref{singtriplo}, ii),
$D'$ contains $2\si$. Let $F\in |F|$ be general: by lemma
\ref{claim7} $g$ maps
$F$ $3$--to--$1$ onto a ruling $R$ of $\f$. Since $F$ is
smooth, by lemma
\ref{singtriplo}
$R$ meets $D'-2\si$ outside $\si$ and by the Hurwitz
formula $(D'-2\si)R=8$. So
$D'-2\si$ is a divisor linearly equivalent to $8\so+aR$,
with $a\ge 0$. We conclude that $D=D'$.

Case ii):

\noindent One proceeds exactly as in case i), setting in
this case
$L=2\si+5R$ and
$M=3\si+7R$. The dimensions of all the linear systems
involved have already been computed in lemma \ref{claim8},
ii).\qed

From the analysis of the linear system $|3C|$ carried out
in section \ref{secg3} it follows that, unless $\k^2=9$ or
$\k^2=6$ and $P\notin Z$, the morphism
$g:\hS\to
\f$ is certainly not finite, and thus the surface $X$ is
not smooth. It is possible to give a precise description
of the curves contracted by $g$, and thus of the singularities
of $X$ in each case, however, since we will not need this, we
will just give sufficient conditions on the singularities
of a cover $\pi:X\to \f$ in order that it arises in the
Stein factorization of $g:\hS\to\f$ as above. In the next
section we will give examples of such singular covers,
thus showing that the surfaces that we have described do
indeed exist and all possible values of the invariants
occur.
\begin{prop}\label{exampleg3}
 Let $\pi:X\to \f$ be a triple cover with trace zero
module $\E$, let $D$ be the branch locus of $\pi$, and let
$h:\hS\to X$ be the minimal desingularization of $X$:
\medskip

\noindent 1) if $\E=
\Of(-2\si-4R)\oplus
\Of(-3\si-8R)$ and  $X$ is nonsingular over $\si$, then the set--theoretic
inverse image of
$\si$ in $\hS$ is an exceptional curve $E$. If  we denote
by $S$ the surface obtained by blowing down
$E$, then the following are true:

1.i) if $X$ is smooth,  then $\k^2=9$, $p_g(S)=5$, and $S$
satisfies assumption
\ref{g3};

1.ii) if  the only singular point of $X$  is a point $x_0$
such that $\pi$ is not totally ramified at $x_0$ and
$y_0=\pi(x_0)$  is a
$(3,3)$ point of $D$,
 then
$\k^2=8$,
$p_g(S)=4$ and
$S$ satisfies assumption
\ref{g3};

1.iii) if  the only singular point of $X$  is a point
$x_0$ such that $\pi$ is not totally ramified at $x_0$ and
$y_0=\pi(x_0)$  is  an ordinary quadruple point of $D$,
then $\k^2=7$, $p_g(S)=4$ and $S$ satisfies assumption
\ref{g3};

1.iv) if $\pi$ is a Galois cover, then
$D= 2\si+2D_0$;  if $D_0$ has an ordinary triple point
$y_0$ and is smooth elsewhere, then $\k^2=6$, $p_g(S)=4$
and $S$ satisfies assumption
\ref{g3} with $P\in Z$ (see lemma \ref{claim5});
\medskip

\noindent 2) if $\E=
\Of(-2\si-5R)\oplus
\Of(-3\si-7R)$ and  $X$ is nonsingular over $\si$, then
the inverse image of
$\si$ in $\hS$ is an exceptional curve $E$, and the
surface $S$ obtained by blowing down $E$ has invariants
$\k^2=6$, $p_g(S)=4$ and satisfies assumption
\ref{g3} with $P\notin Z$ (see lemma \ref{claim5}).
\end{prop}
\proof
  As we have already  seen in the proof of thm.
\ref{Eg3},  in case 1) $X$ is  defined inside the vector
bundle $\E^{\vee}$ by equations
\ref{equazioni}, where

\noindent $a\in H^0(\f, \Of(2\si+4R))$,
 $b\in H^0(\f, \Of(\si))$,

\noindent $c\in H^0(\f, \Of(4\si+12R))$,
 $d\in H^0(\f, \Of(3\si+8R))$;

\noindent so $a$, $b$,
$d$ all vanish on $\si$, and thus $\pi$ is totally
ramified over $\si$ by corollary
$4.6$ of
\cite{mi}. Analogously, in case 2) we have

\noindent$a\in H^0(\f, \Of(2\si+5R))$,
 $b\in H^0(\f,
\Of(\si+3R))$,

\noindent $c\in H^0(\f, \Of(4\si+9R))$,
 $d\in H^0(\f, \Of(3\si+7R))$

 \noindent and so $a$, $c$,
$d$ all vanish on $\si$ and  $\pi$ is totally ramified
over $\si$ also in this case. Thus, both in case 1) and in
case 2), if $X$ is smooth then $h^*\si=3E$, where
$E$ is a smooth  rational curve satisfying
$(3E)^2=3\si^2$, i. e., $E^2=-1$.

Denote by $\hC$ the  pull--back to $\hS$  of the pencil
$|R|$ and by $|C|$ the image of $|\hC|$ in $S$: under our
assumptions, the general $\hC$ (and thus also the general
$C$) is smooth and has genus
$3$  by the Hurwitz formula, since $DR=10$. The
restriction of $\pi\circ h$ maps a smooth
$\hC$
$3$--$1$ onto a smooth rational curve, and it is easy to
check that this implies that
$\hC$ (and thus also
$C$) is not hyperelliptic.
 Moreover the point $P\in S$ to which $E$ is contracted is
the only base point of $|C|$ and it is simple. So $C^2=1$
and, by the adjunction formula,
$\k C=3$.

Next we compute the singularities of $X$ and the
invariants of $S$ in the various cases. In case 1.i),  $X$
is smooth,  and the formulas
of section
$8$ of
\cite{pa}   yield
$p_g(X)=5$ and $K_X^2=8$, and thus $p_g(S)=5$, $\k^2=9$.

In order to compute the invariants of $S$ in  cases 1.ii),
1.iii) and 1.iv), in which $X$ is singular, we describe
$X$ more precisely. Denote by $V$ the open surface
$\f\setminus \si$: $V$ is isomorphic to the total space of
the line bundle $\OO_{\pp^1}(3)$ and, if we denote by
$p:V\to\pp^1$ the projection, then $p_*\OO_V=\oplus _{k\ge
0}\OO_{\pp^1}(-3k)$ and the group
$Pic(V)=p^*Pic(\pp^1)$ is generated by
$L=p^*\OO_{\pp^1}(1)$.
 For
$y\in V$, one has
$b(y)\ne 0$ in equations \ref{equazioni}; thus one can
eliminate $w$ and obtain a relation of the form
$z^3+rz+s=0$, where $r\in H^0(V, L^8)$, $s\in H^0(V,
L^{12})$. Denote again by $L^4$ the total space of the
line bundle $L^4$ on $V$ and by
$q:L^4\to V$ the projection map:  then
$X\setminus \pi\inv(\si)$ is isomorphic to the hypersurface
$X_0=\{z^3+rz+s=0\}\subset L^4$,  where $z$ represents the
tautological section
 $q^*L^4$. We denote again by $\pi:X_0\to V$ the
restriction of $\pi:X\to\f$; notice that
$\pi_*\OO_{X_0}=\OO_V\oplus L^{-4}\oplus L^{-8}$.  The
surface
$S\setminus P$ is   the minimal desingularization of $X_0$
and we denote by
$\eta:S\setminus P\to X_0$ the resolution map.
 By the adjunction formula, the dualizing sheaf
  of $X_0$  is $\omega_{X_0}=\omega_{L^4}\otimes
q^*L^{12}\mid_{X_0}=q^*L^3\mid_{X_0}=\pi^*L^3$; so
$\eta^*\omega_{X_0}=\OO_{S\setminus P}(3C)$ and,
 using the projection formula twice,
$h^0(X_0, \omega_{X_0})=h^0(V, \pi_*\pi^*L^3)= h^0(V,
L^3)=h^0(\pp^1,\OO_{\pp^1}(3))+h^0(\pp^1, \OO_{\pp^1})=5$.

By assumption, the surface
$X_0$ is singular only at the point
$x_0$ such that
$\pi(x_0)=y_0$. The singularity $(X_0, x_0)$ is
analytically isomorphic to the following:

case 1.ii):  $u^2=f(x, y)$, where
$f(x,y)=0$ is the equation of a plane curve having a
$(3,3)$ point at the origin;

case 1.iii): $u^2=f(x, y)$, where $f(x,y)=0$ is the
equation of a plane curve having an ordinary  quadruple
point at the origin;

case 1.iv): $u^3=f(x,y)$, where $f(x,y)=0$ is the equation
of a plane curve having an ordinary  triple  point at the
origin.

In all three cases we have an elliptic Gorenstein
singularity, whose minimal resolution is well known: the
exceptional divisor is an elliptic curve $\Gamma$, such
that $\Gamma^2=-1$ in case 1.i), $\Gamma^2=-2$ in case
1.ii) and
$\Gamma^2=-3$ in case 1.iv).  A section of $\omega_{X_0}$
pulls--back to a regular
$2$-form on $S\setminus P$ iff it vanishes at $x_0$, and
$\omega_{S\setminus P}=\eta^*\omega_{X_0}(-\Gamma)$.
 Since every regular differential form on $S\setminus P$
extends to a global form on $S$, it follows that
$\omega_S=\OO_S(3C-\Gamma)$ and
$p_g(S)=h^0(X_0,
\omega_{X_0})-1=4$.So one has $\k^2=9+\Gamma^2$, namely $\k^2=8$ in case 1.i),
$\k^2=7$ in case 1.ii) and $\k^2=6$ in case 1.iv).

In case 2), $X$ is smooth and   the formulas of section
$8$ of
\cite{pa} yield  $p_g(X)=4$, $K_X^2=5$, and thus
$p_g(S)=4$, $\k^2=6$.

Remark that, since we have $\k^2>0$ and $p_g(S)>0$ in all
cases, the surface  $S$ is of general type.

To show that the surface $S$ satisfies assumption \ref{g3}
in all cases, we only need to show that
$S$ is  minimal. So assume  that this is not the case, and
let $\Delta$ be a $-1$--curve on $S$. Observe that $C$ is
a nef divisor on $S$ and that the only irreducible curves
that have zero intersection with $C$ are those contracted
by the rational map $S \cdots -\! >\f$. By the above
discussion, no rational curve on
$S$ is contracted,  and so $\Delta C=m>0$. Let $Y$ be the
surface obtained from $S$ by contracting $\Delta$ and let
$D$ be the image of $C$ in $Y$; one has
$D^2=C^2+m^2=1+m^2$, $K_Y D=\k C-m=3-m$. So the index
theorem yields: $2K_Y^2\le D^2 K_Y^2\le (K_YD)^2\le 4$,
which contradicts $K_Y^2=\k^2+1\ge 7$.

Finally, by prop. \ref{Eg3}, i) we have $P\in Z$ in case
1.iv) and $P\notin Z$ in case 2).\qed

 Next we repeat the previous analysis for surfaces as in
assumption \ref{g4} with $\deg\phi_K=3$.
\begin{nota}\label{notaII} Let $S$ be a surface satisfying
assumption \ref{g4} and such that
$\deg\phi_K=3$, let
$\epsi':S'\to S$,
$C'$,
$E'$ as in notation
\ref{blowup}, let $\epsi:\hS\to S'$ be the blow--up of
$S'$ at the point $P_2$,
 (see lemma \ref{basepointg4}), let $E$ be the
exceptional curve of $\epsi$, let $C'$ denote again the
pull--back of $C'$ on $\hS$ and let $F=C'-E$.
\end{nota}
\begin{prop}\label{morfismog4} Let $f:S'\to C_3$ be the
morphism defined in lemma \ref{tripleg4}; using notation
\ref{notaII}, there is a commutative diagram:
\begin{equation}
\begin{array}{rcccl}
\phantom{1} &\hS&\stackrel{\epsi}{\rightarrow} &S'&
\phantom{1} \\
\scriptstyle{g}\!\!\!\!\!\! & \downarrow & \phantom{1} &
\downarrow &
\!\!\!\!\!\!
\scriptstyle{f}
\\
\phantom{1} & \f & \stackrel{\psi}{\rightarrow} & C_3 &
\phantom{1}
\end{array}
\end{equation} where $g:\hS\to\f$ is a finite flat
morphism of degree $3$. The map
$g$ is totally ramified and finite over
$\si$.
\end{prop}
\proof The proof is the same as the proof of prop.
\ref{morfismo}. In this case
$g^*(R+\so)=F+3C'$ and, by lemma \ref{dimg4}, i), the
system $|F+3C'|$ is free of dimension  $7=\dim |R+\so|$.
So $g^*|R+\so|=|F+3C'|$ is free and $g$ is a morphism. The
last part of the statement is proven as in prop.
\ref{morfismo}.\qed

\begin{nota-rem}\label{steing4}If $S$ is a surface as in
assumption
\ref{g4} and such that $\deg\phi_K=3$, $\hat{S}$ is as in
notation \ref{blowup}, we denote by
$\hS\stackrel{h}{\to}X\stackrel{\pi}{\to}\f$  the Stein
factorization of
$g:\hS\to \f$ (see prop. \ref{morfismog4}), so that $X$ is
a normal variety, $h$ is a morphism with connected fibres
and
$\pi:X\to\f$ is a finite  morphism of degree $3$.  Notice
that by lemma
\ref{flatness}, $\pi$ is flat.  So  the theory of
\cite{mi} and \cite{pa} applies and
$\pi_*\OX=\OO_{\f}\oplus \E$, where $\E$ is a locally free
sheaf of rank
$2$, the so-called {\em trace--zero module} of $\pi$.
\end{nota-rem} In the next proposition we determine the
trace--zero module $\E$.
\begin{prop} \label{Eg4} Let $S$ be a surface as in
assumption \ref{g4} with $\deg\phi_K=3$ and let
$\pi:X\to
\f$ and ${\cal E}$ be as in notation--remark
\ref{steing4}.  Then
$\E=\OO_{\f}(-2\si-5R)\oplus \OO_{\f}(-4\si-10R)$.
\end{prop}
\proof The proof is analogous to that of  prop. \ref{Eg3}
and therefore we will only sketch it.

Set $L=\Of(2\si+5R)$; then by lemma \ref{tripleg4}
$g^*L=\OSh(6E+5F)=\OSh(5C'+E)$ and $|g^*L|=E+|5C'|$. So
lemma \ref{dimg4} ii) implies  $\dim|g^*L|=9$ and the
vanishing order on $E$ of a general element  of
$H^0(\hS,g^*L)$ is precisely $1$. A general element $z$ of
$H^0(\hS,g^*L)$ defines a birational map $\psi:\hS\to L$
such that $g$ factors through $\psi$. Arguing as in the
proof of prop. \ref{Eg3} and using  lemma \ref{dimg4}, one
can show that $H^0(\hS,g^*L^3)=g^*H^0(\f,L)\oplus
zH^0(\f,L^2)\oplus z^2H^0(\f,L^3)$. So there exist $a_i\in
g^*H^0(\f,L^i)$, $i=1,2,3$ such that
$z^3+a_1z^2+a_2z+a_3=0$; up to replacing $z$ by
$z+\frac{a_1}{3}$, one may assume that $z^3+rz+s=0$, where
$r\in  g^*H^0(\f,L^2)$ and $s\in  g^*H^0(\f,L^3)$.
Reasoning again as in the proof of prop. \ref{Eg4} one
shows that the image of
$\psi$ is a normal surface defined inside $L$ by the
equation $z^3+rz+s=0$ and therefore it is equal to $X$ and
$E= L\inv\oplus L^{-2}$. \qed

We close this section by stating the analogue of prop.
\ref{exampleg3}:
\begin{prop}\label{exampleg4} Let $\pi:X\to \f$ be a
triple cover with trace zero module
$\E=\OO_{\f}(-2\si-5R)\oplus \OO_{\f}(-4\si-10R)$, let $D$
be the branch locus of
$\pi$ and let
$h:\hS\to X$ be the minimal desingularization of $X$.

If $X$ is nonsingular over $\si$, then  the inverse image
of
$\si$ in $\hS$ is an exceptional curve $E$; denote by $S'$
the surface obtained by blowing down
$E$ and by $P_2$ the image point of $E$. If moreover $X$
has only a singular point
$x_0$ such that $\pi$ is not totally ramified at $x_0$ and
$y_0=\pi(x_0)$
 is an  ordinary octuple point of $D$, then $S'$ contains
an exceptional curve
$E'$ such that $P_2\in E'$, and the surface $S$ obtained
from $S'$ by contracting
$E'$ has invariants
 $\k^2=8$, $p_g(S)=4$, satisfies assumption
\ref{g4} and the canonical map of $S$ has degree $3$.
\end{prop}

\proof The first statement can be proven exactly as in
prop.
\ref{exampleg3}.

The linear system $|\hC|=h^*\pi^*|R|$ is a base point
free  pencil of curves of genus
$4$, such that the generic element is smooth; moreover, as
it was the case in prop.
\ref{exampleg3}, the smooth curves of the pencil are not
hyperelliptic, since they admit a morphism of degree
$3$ onto a rational curve. The image $|C'|$ of $|\hC|$ on
$S'$ has the point
$P_2$ to which
$E$ is contracted as a simple base point, and thus
$C'^2=1$, $K_{S'}C'=5$.

Under our assumptions on $\E$, the term $b$ in equations
\ref{equazioni} is a constant; moreover
$b\ne 0$, since otherwise  $D_0$ would not be reduced by
prop. $4.5$ of
\cite{mi}. So, as in the proof of prop. \ref{exampleg3},
one can eliminate $w$ from equations
\ref{equazioni} and obtain a relation of the form
$z^3+rz+s=0$, where $r\in H^0(\f, \Of(4\si+10R))$ and
$s\in H^0(\f, \Of(6\si+15R))$. Thus if we denote by $L$
the total space of the line bundle $\Of(2\si+5R)$ on $\f$
and by $q:L\to
\f$ the projection, then
 $X$ is isomorphic to the hypersurface
$\{z^3+rz+s=0\}\subset L$ , where
$z$ represents the tautological section of $q^*L$. By the
adjunction formula, the dualizing sheaf of $X$ is
$\omega_X=\omega_{L}\otimes q^*
L^3\mid_X=q^*(\omega_{\f}\otimes
L^2)|_X=\pi^*\Of(2\si+5R)$;  using the projection formula,
one computes
$h^0(X,\omega_X)=h^0(\f, \Of)+h^0(\f,
\Of(2\si+5R))=10$. Denoting the resolution map by
$\eta:\hS\to X$, one has
$\eta^*\omega_X=\OO_{\hS}(6E+5\hC)$.

The singularity $(X,x_0)$ is analytically isomorphic to
the hypersurface singularity $u^2+f(x,y)=0$, where
$f(x,y)=0$ is the equation of a plane curve with an
octuple point at the origin. The minimal resolution of
this singularity can be computed by blowing--up the
$x,y$--plane at the origin and then taking pull--back and
normalization. The exceptional divisor is a smooth
hyperelliptic curve $\Gamma$ of genus $3$ such that
$\Gamma^2=-2$ and one has
$\omega_{\hS}=\eta^*\omega_X(-3\Gamma)$, so that
$K_{\hS}^2=6$. The condition that the pull-back of a
section $\sigma\in H^0(X,\omega_X)$ be a regular form on
$\hS$ is expressed in local coordinates by requiring
$\frac{\partial^{i+j}\sigma}{\partial x^i y^j}(x_0)=$ for
$0\le i+j\le 2$; this amounts to $6$ linear conditions and
thus
$p_g(\hS)\ge h^0(X,\omega_X)-6=4$. Denote by $R_0$ the
ruling of $\f$ that contains $y_0$ and write
$\eta^*(\pi^*R_0)=\Gamma +Z$; since $\Gamma+Z\simeq \hC$,
one has $\Gamma Z=2$,
$Z^2=-2$, $K_{\hS}Z=0$. Moreover,  $Z$ is a  finite triple
cover of $R_0$ totally ramified over the intersection
point of $\si$ and $R_0$, and the total ramification point
is a smooth point of $Z$: it follows easily that $Z$ is
irreducible and therefore smooth, since $p_a(Z)=0$.  Since
$Z$ meets $E$ transversally, the image $E'$ of $Z$ on $S'$
is an exceptional curve containing
$P_2$, and therefore $E'$ can be contracted to a smooth
point $P$, yielding a surface
$S$ with
$p_g(S)=p_g(\hS)\ge 4$ and $\k^2=8$. The image of the
pencil $|C'|$ is a pencil
$|C|$ on $S$, satisfying $C^2=2$, $\k C=4$. The surface
$S$ is of general type, since $p_g(S)>0$ and $\k^2>0$. In
order to show that $S$ is minimal, one can argue exactly
as in the proof of thm. \ref{exampleg3}. This completes
the proof that $S$ satisfies assumption \ref{g4}. Thus, by
lemma \ref{1}, $p_g(S)=4$. Again by lemma \ref{1}, in
order to show that $\deg\phi_K=3$ it is enough to exclude
that $\deg\phi_K=4$. So assume by contradiction that
$\deg\phi_K=4$, and denote by
$\Sigma$ the image of $\phi_K$ and by $d$ the degree of
$\Sigma$: we have
$8=\k^2\ge 4d\ge 4(p_g(S)-2)=8$, and therefore $d=2$ and
$\phi_K$ is a morphism. Recall that by lemma \ref{1}
$\k=2C$. Consider a general $C$: the system
$|3P|$ is free of dimension
$1$ by assumption,  and the system $|4P|$ is  also free,
since it contains the free system $|\k|_C$. It follows
that $h^0(C,\OO_C(4P))=3$. On the other hand, we have
$K_C=6P$ by the adjunction formula,  and therefore the
Riemann-Roch theorem gives:
$h^0(C,\OO_C(4P))=1+h^0(C,\OO_C(2P)=2$. So we have reached
a contradiction.
\qed

\section{The examples}
\setcounter{defn}{0}
\setcounter{equation}{0}

In this section we  give explicit examples of surfaces satisfying
assumption
\ref{g3} and taking all the possible values of the invariants (see thm.
\ref{I}), and of a surface satisfying assumption \ref{g4} (in this case the
only possibility for the invariants is $\k^2=8$, $p_g(S)=4$, by thm. \ref{II}).
We do this by showing that the assumptions in prop.
\ref{exampleg3} and
\ref{exampleg4} can actually be verified.
 The cases in
prop.
\ref{1i},
\ref{2} and \ref{1iv} are shown to exist by Bertini type arguments.

\begin{prop}\label{1i}
There exists a surface $S$ satisfying assumption \ref{g3} with invariants
$\k^2=9$, $p_g(S)=5$.
\end{prop}
\proof
It is enough to show the existence of a smooth
triple cover $\pi:X\to\f$ as in  prop. \ref{exampleg3}, 1.i).
As we have already remarked several times, such a cover is determined by the
choice
of $a\in H^0(\f,\Of(2\si+4R))$, $b\in H^0(\f,\Of(\si))$,
$c\in H^0(\f,\Of(4\si+12R))$, $d\in H^0(\f,\Of(3\si+8R))$. If one  takes $a=0$,
$d=0$, $b\ne 0$ and $c$ such that $c=0$ is a smooth divisor, then it is easy to
check using lemma \ref{singtriplo} and equations \ref{equazioni} that the
corresponding cover is smooth. Notice that by Bertini's theorem  it is
possible to
find
$c$ as required since the linear system $|4\so|$ is base point free.\qed

\begin{prop}\label{1iv} There exists a surface $S$ satisfying assumption
\ref{g3}
with invariants $\k^2=6$,
$p_g(S)=4$ and $P\in Z$ (see lemma \ref{claim5}).
\end{prop}
\proof By prop. \ref{exampleg3}, 1.iv), we have to show that there exists a
Galois triple
cover $\pi:X\to \f$ branched over $\si+D_0$, where $D_0\in |4\so|$ has an
ordinary
triple point $y_0\notin \si$ and is smooth elsewhere. This corresponds to
taking,
in equations \ref{equazioni}, $a=d=0$, $b\in H^0(\f,\Of(\si))\setminus \{0\}$,
and
$c\in H^0(\f,\Of(4\so))$ such that the divisor $D_0=\{c=0\}$ is as above.
If the point $y_0$ is fixed, then the linear system $W$ of divisors in $|4\so|$
having a triple point at $y_0$ has no extra base point and the generic
element of
$W$ has an ordinary triple point at $y_0$: this can be seen easily by
considering
reducible elements in $W$ given by the sum of three  curves of $|\so|
$ passing
through $y_0$ and of a fourth curve of $|\so|$ not passing through $y_0$.
So by Bertini's theorem the generic divisor of $W$ has an ordinary triple point
at $y_0$ and is smooth elsewhere.\qed

\begin{prop}\label{2} There exists a surface $S$ satisfying assumption \ref{g3}
with invariants
$\k^2=6$, $p_g(S)=4$ and $P\notin Z$ (see lemma \ref{claim5}).
\end{prop}
\proof By prop. \ref{exampleg3}, 2), we have to show that there exists a triple
cover
$\pi:X\to \f$ with $X$ smooth and trace zero module
$\E=\Of(-2\si-5R)\oplus\Of(-3\si-7R)$. Denote by $\pp$ the $\pp^1$--bundle
$\hbox{Proj}(\E^{\vee})$ on $\f$, by $T$ the tautological hyperplane section on
$\pp$ and by $p:\pp\to \f$ the projection. Let $X\subset \pp$ be a smooth
divisor
linearly equivalent to $3T+p^*\det\E$, and assume moreover that $X$ contains no
fibre of the map $p$: then the map $\pi:X\to \f$, induced by restricting
$p$, is quasi--finite,  and therefore flat by lemma \ref{flatness}. Thus
$\pi:X\to\f$ is a triple cover and, by prop. $8.1$ of \cite{mi}, its trace zero
module is equal to
$\E$. So, in order to prove the claim, it is enough to show that the linear
system
$|3T+p^*\det\E|$ contains such a divisor $X$.

Recall that
$p_*\OO_{\pp}(3T+\det\E)=\Of(\si+3R)\oplus
\Of(2\si+5R)\oplus\Of(3\si+7R)\oplus\Of(4\si+9R)$: since the general section of
this vector bundle  vanishes nowhere on $\f$, the general element of
$|3T+p^*\det\E|$ contains no fibre of $p$. Moreover, it is not difficult to
check that, as a set,  the base locus $B$ of
$|3T+p^*\det\E|$ is the intersection of $p\inv(\si)$ with  the section of $\pp$
corresponding to the quotient map $\E^{\vee}\to\Of(2\si+5R)$
 and that the general divisor of $|3T+p^*\det\E|$  is smooth at every
point of $B$.
 It follows from Bertini's theorem that the general divisor in
$|3T+p^*\det\E|$ is smooth.\qed

\bigskip

The remaining examples are constructed explicitly,  using the symbolic
computation
program Axiom to show that the  assumptions of  propositions
\ref{exampleg3} and
\ref{exampleg4} are satisfied.
 We start by setting some notation. We denote by $(t_0:t_1)$
homogeneous coordinates on $\pp^1$, by $U_i$ the open subset $\{t_i\ne
0\}\subset
\pp^1$,
$i=0,1$, by $t=\frac{t_1}{t_0}$ the affine coordinate on $U_0$ and by
$s=\frac{t_0}{t_1}$ the affine coordinate on $U_1$.
The open surface $V=\f\setminus\si$ is the union of two affine open subsets
$V_i=U_i\times \C$, $i=0,1$, with coordinates $(t,u)$ and $(s,v)$
respectively, related   by $t=\frac{1}{s}$, $u=s^3v$ on $V_0\cap V_1$.

\begin{prop}\label{1ii}
There exists a surface $S$ satisfying assumption \ref{g3} with invariants
$\k^2=8$,
$p_g(S)=4$.
\end{prop}
\proof
It is enough to show the existence of a triple cover
$\pi:X\to\f$ as in prop. \ref{exampleg3}, 1.ii). Such a cover is determined
by the
choice of of $a\in H^0(\f,\Of(2\si+4R))$, $b\in H^0(\f,\Of(\si))$,
$c\in H^0(\f,\Of(4\si+12R))$, $d\in H^0(\f,\Of(3\si+8R))$. We take $a=0$,
$b\ne 0$;
using lemma \ref{singtriplo} and lemma $4.5$ of \cite{mi}, it is easy to
check that
if $c$ does not vanish on $\si$ then $X$ is nonsingular over $\si$. Therefore we
only need to study the singularities of $X$ over $V=\f\setminus\si$.  We may
assume that $b=1$ on $V$ and eliminate $w$ in equations \ref{equazioni}.
Denote by
$q:V\to\pp^1$ the projection and by $L$ the line bundle $q^*\OO_{\pp^1}(4)$ on
$V$: then $\pi\inv(V)$ is isomorphic to the hypersurface $\{z^3+3dz-c=0\}\subset
L$, where $z$ represents the tautological section of $L$. On the open set
$V_0$ we
can write: $d=\sum_{i=0}^2d_iu^{2-i}$, where the $d_i$'s are polynomials in
$t$ of
degree $2+3i$, and $c=\sum_{i=0}^4c_iu^{4-i}$, where the $c_i$'s are
polynomials in
$t$ of degree $3i$. The condition that $c$ does not vanish on $\si$
corresponds to
the condition $c_0\ne 0$. Moreover, we assume that the image of the
singular point
$x_0$ is the point $y_0=(0,0)\in V_0$ (this can always be achieved by means
of an
automorphism of
$\f$). The computations of section \ref{axiom} show that the required
example can
be obtained, for instance, by choosing
$d=(1+t^2)u^2+2u-t^8+t^7-t^6-1$ and $c=-2u^4-2u^3+6u-2$.
Although we
will not give the computations here, we have used Axiom also to determine the
coefficients of $d$ and $c$ in such a way that the branch locus $D$ of the cover
has  a $(3,3)$--point at $y_0$.\qed

\begin{prop}\label{1iii} There exists a surface $S$ satisfying assumption
\ref{g3}
with invariants $\k^2=7$,
$p_g(S)=4$.
\end{prop}
\proof
It is enough to show the existence of a triple cover
$\pi:X\to\f$ as in prop. \ref{exampleg3}, 1.iii). We argue as in the proof prop.
\ref{1ii} and, using the same notation, we take here: $a=0$, $b\ne 0$,
$c=-2u^4-20u^3-6u^2+12u+2t^{12}-2t^{11}-2$, $d=2u^2+4u+t^4-1$.
Also in this
case,
we have used Axiom  to determine the coefficients of $d$ and $c$ in such a way
that the branch locus $D$ of the cover has  an ordinary quadruple point at
$y_0$.\qed

\begin{prop}\label{due}
 There exists a surface $S$ satisfying assumption \ref{g4}.
As predicted  by thm. \ref{II}, the invariants of $S$ are $\k^2=8$, $p_g(S)=4$.
\end{prop}
\proof
It is enough to show the existence of a triple cover
$\pi:X\to\f$ as in prop. \ref{exampleg4}; using a notation consistent with the
one in the proof of prop. \ref{exampleg4}, we wish to find $r\in H^0(\f,
\Of(4\si+10R)$ and $s\in H^0(\f,\Of(6\si+15R)$ such that the hypersurface
$X=\{z^3+rz+s=0\}\subset \Of(2\si+5R)$ satisfies the assumptions of prop.
\ref{exampleg4}. Notice that $X$ is smooth over $\si$ iff $s$
vanishes on $\si$ of order $1$. So, as in the previous cases, we only  study
the restriction of $X$ to $V$ and assume that $y_0=(0,0)\in V_0$. Using the
coordinates of $V_0$, we write $r=\sum_{i=0}^3r_iu^{3-i}$, where $r_i$ is a
polynomial in $t$ of degree $3i+1$ and $s=\sum_{i=0}^5s_iu^{5-i}$, where $s_i$
is a polynomial in $t$ of degree $3i$. By the computations of section
\ref{axiom}, we may take:
$r=36u^3-45u^2+18u-3+3t^{10}-3t^9+3t^8$ and
$s=-27u^5+135u^4-144u^3+72u^2-18u+2$.
 Also in this case, we have used Axiom  to
determine the coefficients of $d$ and
$c$ in such a way that the branch locus $D$ of the cover has  an ordinary
octuple point at
$y_0$.\qed

\section{Appendix: computations with Axiom}\label{axiom}
\setcounter{defn}{0}
\setcounter{equation}{0}

This section contains the computations with Axiom that are needed in the proofs
of propositions
\ref{1ii}, \ref{1iii} and
\ref{due}. We give a slightly edited version of the Axiom session for prop.
\ref{due}, which shows the existence of   the ``new'' example, and only the
input sequences for prop.
\ref{1ii} and
\ref{1iii}, which are very similar to the first one. The notation is consistent
with the one defined in the previous section, the only difference is that in the
first example the affine coordinate at infinity on $\pp^1$ is denoted by $y$
instead of $s$.

{\underline Prop. \ref{due}}:

\noindent We start by checking that the only singular point of $X$ over
$U_0$ is the point $z=1, u=t=0$

\begin{verbatim}
initial (60) -> r:P := 36*u^3-45*u^2+18*u-3+3*t^10-3*t^9+3*t^8
\end{verbatim}

\noindent(60)\quad  $36u^3 - 45u^2  +18u + 3t^{10}   - 3t^9  + 3t^8  - 3$

\begin{verbatim}
                                      Type: Polynomial Fraction Integer
initial (61) -> s:P :=-27*u^5+135*u^4-144*u^3+72*u^2-18*u+2
\end{verbatim}

\noindent(61)\quad  $- 27u^5  + 135u^4
-144u^3+72u^2-18u+2$

\begin{verbatim}
                                      Type: Polynomial Fraction Integer
initial (62) -> f:P:=z^3+r*z+s
\end{verbatim}

\noindent(62)\quad $ z^3  + (36u^3  - 45u^2 + 18u + 3t^{10}  - 3t^9  + 3t^8
-3)z - 27u^5  +
135u4  - 144u^3    + 72u^2  - 18u + 2$
\begin{verbatim}

                                      Type: Polynomial Fraction Integer
initial (63) -> fz:P:=differentiate(f,z)
\end{verbatim}

\noindent(63)\quad  $3z^2  + 36u^3  - 45u^2  + 18u + 3t^{10}   - 3t^9  +
3t^8  - 3$

\begin{verbatim}
                                      Type: Polynomial Fraction Integer
initial (64) -> ft:P:=differentiate(f,t)
\end{verbatim}

\noindent(64)\quad $(30t^9  - 27t^8  + 24t^7 )z$

\begin{verbatim}
                                      Type: Polynomial Fraction Integer
initial (65) -> fu:P:=differentiate(f,u)
\end{verbatim}

\noindent(65)\quad  $(108u^2  - 90u + 18)z - 135u^4  + 540u^5  - 432u^2  +
144u - 18$

\begin{verbatim}
                                      Type: Polynomial Fraction Integer
initial (66) -> sing:List P :=[f,fz,ft,fu]
\end{verbatim}

\noindent(66)\quad $[ z^3  + (36u^3  - 45u^2 + 18u + 3t^{10}  - 3t^9  +
3t^8 -3)z - 27u^5  +
135u4  - 144u^3    + 72u^2  - 18u + 2,3z^2  + 36u^3  - 45u^2  + 18u + 3t^{10}
- 3t^9  + 3t^8  - 3,(30t^9  - 27t^8  + 24t^7 )z,(108u^2  - 90u + 18)z - 135u^4
+ 540u^5  - 432u^2  + 144u - 18]$

\begin{verbatim}
                                  Type: List Polynomial Fraction Integer

initial (67) -> solve sing
+++ Garbage collection 27  (internal list2*) after 177.89+65.34 seconds
At gc end about 5.5Mbytes of 14.0 (39.3%) of heap is in use
(67)  [[z= 1,u=0,t= 0]]
                     Type: List List Equation Fraction Polynomial Integer
\end{verbatim}

\noindent Next we compute the discriminant of $f$ with respect to the variable
$z$, which is  an equation for the branch locus of the cover
$X|_{U_0}\to U_0$. It is easy to check that it has an ordinary octuple point at
$t=u=0$.

\begin{verbatim}
initial (68) -> D:P :=4*r^3+27*s^2
\end{verbatim}

\noindent(68)\quad $19683u^{10}   - 10206u^9  + 2187u^8  + (46656t^{10}   -
46656t^9  +
46656t^8 )u^6    +$

$(-116640t^{10}+116640t^9-116640t^8
)u^5+(119556t^{10}-119556t^9+119556t^8)u^4+$

$+(3888t^{20}-7776t^{19}+11664t^{18}-7776t^{17}+3888t^{16}-66096t^{10}+66096t^9
 - 66096t^3)u^3+$

$+(-4860t^{20}+ 9720t^{19}-14580t^{18}+
9720t^{17}-4860t^{16}+21384t^{10}-21384t^9+21384t^8)u^2+$

$+(1944t^{20}-3888t^{19}+5832t^{18}-3888t^{17}+1944t^{16}-3888t^{10}  +
3888t^9-3888t^8)u+$

$+108t^{30}-324t^{29}+648t^{28}-756t^{27}+648t^{26}-324t^{25}+108t^{24}-
324t^20+$

$+648t^{19}-972t^{18}+648t^{17}-324t^{16}+324t^{10}-324t^9+324t^8$

\begin{verbatim}
                                      Type: Polynomial Fraction Integer
\end{verbatim}

\noindent Finally we show that $X$ is smooth above the line $\{y=0\}\subset
U_1$, where $y=1/t$ is the affine coordinate at infinity on $\pp^1$. We denote
by $f_0=z^3+r_0z+s_0$ the equation of $X$ over the open set $U_1$.

\begin{verbatim}
initial (76) -> r0:P:=36*y*u^3-45*y^4*u^2+18*y^7*u-3*y^10+3-3*y+3*y^2
\end{verbatim}

\noindent(76)\quad  $- 3y^{10}   + 18u y^7  - 45u^2 y^4  + 3y^2  + (36u^3
-3)y+3$

\begin{verbatim}
                                       Type: Polynomial Fraction Integer
initial (77)->s0:P:=-27*u^5+135*y^3*u^4-144*y^6*u^3+72*y^9*u^2-18*y^12*u+2*y
\end{verbatim}

\noindent(77)\quad  $- 18uy^{12}+72u^2 y^9-144u^3 y^6+135u^4 y^3  + 2y - 27u^5$

\begin{verbatim}
                                      Type: Polynomial Fraction Integer
initial (78) -> f0:P:=z^3+r0*z+s0
\end{verbatim}

\noindent(78)\quad $z^3+(-3y^{10}+18u y^7- 45u^2y^4  + 3y^2  + (36u^3  -
3)y + 3)z - 18u
y^{12}   + 72u^2 y^9
 - 144u^3 y^6  + 135u^4 y^3  + 2y - 27u^5$

\begin{verbatim}
                                      Type: Polynomial Fraction Integer
initial (79) -> f0u:P :=differentiate(f0,u)
\end{verbatim}

\noindent(79)\quad  $(18y^7  - 90u y^4  + 108u^2 y)z - 18y^{12}  + 144u y^9
- 432u^2 y^6  +
540u^3 y^3  - 135u^4$

\begin{verbatim}
                                      Type: Polynomial Fraction Integer
initial (80) -> f0y:P :=differentiate(f0,y)
\end{verbatim}

\noindent(80)\quad $(- 30y^9  + 126u y^6  - 180u^2 y^3  + 6y + 36u^3  - 3)z
- 216u y^{11}   +
648u^2 y^8    +
     - 864u^3 y^5  + 405u^4 y^2 +2$

\begin{verbatim}
                                      Type: Polynomial Fraction Integer
initial (81) -> f0z:P :=differentiate(f0,z)
\end{verbatim}

\noindent(81)\quad  $3z^2  - 3y^{10}   + 18uy^7  - 45u^2 y^4  + 3y^2  +
(36u^3  - 3)y + 3$

\begin{verbatim}
                                      Type: Polynomial Fraction Integer
initial (82) -> sing0:List P:=[f0,f0z,f0y,f0u,y]
\end{verbatim}

\noindent(82)\quad    $[  z^3+(-3y^{10}+18u y^7- 45u^2y^4  + 3y^2  + (36u^3
- 3)y + 3)z - 18u
y^{12}   + 72u^2 y^9
 - 144u^3 y^6  + 135u^4 y^3  + 2y - 27u^5, 3z^2  - 3y^{10}   + 18uy^7  - 45u^2
y^4  + 3y^2  + (36u^3  - 3)y + 3,(- 30y^9  + 126u y^6  - 180u^2 y^3  + 6y +
36u^3  - 3)z - 216u y^{11}   + 648u^2 y^8    +
     - 864u^3 y^5  + 405u^4 y^2 +2,(18y^7  - 90u y^4  + 108u^2 y)z - 18y^{12}  +
144u y^9  - 432u^2 y^6  + 540u^3 y^3  - 135u^4, y]$

\begin{verbatim}
                                  Type: List Polynomial Fraction Integer
initial (83) -> solve %
(83)  [[]]                        Type: List List Equation Fraction
Polynomial Integer
\end{verbatim}

{\underline Prop. \ref{1ii}}:
\begin{verbatim}
initial (44) -> c:P :=-2*u^4-2*u^3+6*u-2
initial (45) -> d:P := (1+t^2)*u^2+2*u-t^8+t^7-t^6-1
initial (46) -> f:P := z^3+3*d*z-c
initial (47) -> fz:=differentiate(f,z)
initial (48) -> fu:=differentiate(f,u)
initial (49) -> ft:=differentiate(f,t)
initial (50) -> sing:List P :=[f,fu,ft,fz]
initial (51) -> solve sing
initial (52) -> D:P := 4*d^3 +c^2
initial (53) -> d0:P:= (s^2+1)*u^2+2*s^5*u-1+s-s^2-s^8
initial (54) -> c0:P := -2*u^4-2*s^3*u^3+6*s^9*u-2*s^12
initial (55) -> f0:P := z^3+3*d0*z-c0
initial (56) -> f0z:P :=differentiate(f0,z)
initial (57) -> f0s:P :=differentiate(f0,s)
initial (58) -> f0u:P :=differentiate(f0,u)
initial (59) -> sing0: List P :=[f0,f0z,f0u,f0s,s]
initial (60) -> solve sing0
\end{verbatim}

{\underline Prop. \ref{1iii}}:
\begin{verbatim}
initial (26) -> d:P := 2*u^2+4*u+t^4-1
initial (27) -> c:P := -2*u^4-20*u^3-6*u^2+12*u+2*t^12-t^11-2
initial (28) -> f:P :=z^3+3*d*z-c          3
initial (29) -> fz:=differentiate(f,z)
initial (30) -> fu:=differentiate(f,u)
initial (31) -> ft:=differentiate(f,t)
initial (32) -> sing:List P :=[f,fu,ft,fz]
initial (33) -> solve sing
initial (34) -> D:P := 4*d^3 +c^2
initial (35) -> d0:P := 2*s^2*u^2+4*s^5*u+s^4-s^8
initial (36) -> c0:P := -2*u^4-20*s^3*u^3-6*s^6*u^2+12*s^9*u+2-2*s-2*s^12
initial (37) -> f0:P := z^3+3*d0*z-c0
initial (38) -> f0s:P :=differentiate(f0,s)
initial (39) -> f0z:P :=differentiate(f0,z)
initial (40) -> f0u:P :=differentiate(f0,u)
initial (41) -> sing0 :List P :=[f0,f0z,f0u,f0s,s]
initial (42) -> solve sing
\end{verbatim}

\bigskip

\begin{tabbing}
1699 Lisboa Codex, PORTUGALxxxxxxxxx\= 56127 Pisa, ITALY \kill
Margarida Mendes Lopes               \> Rita Pardini\\
Centro de Algebra \> Dipartimento di Matematica\\
Universidade de Lisboa \> Universit\a`a di Pisa \\
Av. Prof. Gama Pinto, 2 \> Via Buonarroti 2\\
1699 Lisboa Codex, PORTUGAL \> 56127 Pisa, ITALY\\
mmlopes@lmc.fc.ul.pt \> pardini@dm.unipi.it
\end{tabbing}

 \end{document}